\documentclass{article}

\usepackage[utf8]{inputenc}
\usepackage[T1]{fontenc}
\usepackage[margin=1.5in]{geometry}
\usepackage[hidelinks]{hyperref}
\usepackage[x11names, table]{xcolor}
\usepackage{booktabs}
\usepackage{tikz}
\usepackage{pgfplots}
\pgfplotsset{compat=1.10}
\usepackage{dsfont}
\usepackage{amsmath}
\usepackage{amssymb}
\usepackage{amsfonts}
\usepackage{mathtools}

\usepackage{amsthm}
\newtheorem{theorem}{Theorem}[section]
\newtheorem{definition}[theorem]{Definition}
\newtheorem{proposition}[theorem]{Proposition}

\newcommand{\imag}{\boldsymbol{\imath}}

\definecolor{tab_blue}{HTML}{1F77B4}
\definecolor{tab_orange}{HTML}{FF7F0E}

\newcommand{\bbR}{\ensuremath{\mathbb{R}}}
\newcommand{\bbC}{\ensuremath{\mathbb{C}}}
\newcommand{\Rn}{\ensuremath{\bbR^{n}}}
\newcommand{\Rm}{\ensuremath{\bbR^{m}}}
\newcommand{\Rp}{\ensuremath{\bbR^{p}}}

\newcommand{\Rnn}{\ensuremath{\bbR^{n \times n}}}
\newcommand{\Rnm}{\ensuremath{\bbR^{n \times m}}}
\newcommand{\Rpn}{\ensuremath{\bbR^{p \times n}}}

\newcommand{\Rmp}{\ensuremath{\bbR^{m \times p}}}

\newcommand{\cA}{\mathcal{A}}
\newcommand{\cB}{\mathcal{B}}
\newcommand{\cC}{\mathcal{C}}
\newcommand{\cE}{\mathcal{E}}
\newcommand{\cH}{\mathcal{H}}
\newcommand{\cL}{\mathcal{L}}
\newcommand{\cM}{\mathcal{M}}
\newcommand{\cN}{\mathcal{N}}
\newcommand{\cS}{\mathcal{S}}
\newcommand{\cT}{\mathcal{T}}
\newcommand{\cV}{\mathcal{V}}
\newcommand{\cW}{\mathcal{W}}

\newcommand{\hH}{\hat{H}}

\newcommand{\hx}{\hat{x}}
\newcommand{\hy}{\hat{y}}

\newcommand{\hcA}{\hat{\cA}}
\newcommand{\hcB}{\hat{\cB}}
\newcommand{\hcC}{\hat{\cC}}
\newcommand{\hcE}{\hat{\cE}}
\newcommand{\hcH}{\hat{\cH}}

\newcommand{\one}{\mathds{1}}

\newcommand{\fundef}[3]{{#1}\colon{#2}\to{#3}}

\newcommand{\graph}{\mathsf{G}}
\newcommand{\vset}{\mathsf{V}}
\newcommand{\ver}{\mathsf{v}}
\newcommand{\vnum}{\mathsf{n}}
\newcommand{\VL}{\vset_{\operatorname{L}}}
\newcommand{\VF}{\vset_{\operatorname{F}}}
\newcommand{\eset}{\mathsf{E}}
\newcommand{\edge}{\mathsf{e}}
\newcommand{\wfun}{\mathsf{w}}
\newcommand{\weight}{\mathsf{a}}

\newcommand{\amat}{\mathsf{A}}

\newcommand{\degree}{\delta}
\newcommand{\dmat}{\mathsf{D}}
\newcommand{\lmat}{\mathsf{L}}

\newcommand{\rmat}{\mathsf{R}}
\newcommand{\wmat}{\mathsf{W}}
\newcommand{\mass}{\mathsf{m}}

\newcommand{\mmat}{\mathsf{M}}

\newcommand{\imat}{\mathsf{B}}

\newcommand{\iw}{\mathsf{b}}

\newcommand{\inum}{\mathsf{m}}
\newcommand{\omat}{\mathsf{C}}

\newcommand{\ow}{\mathsf{c}}

\newcommand{\onum}{\mathsf{p}}
\newcommand{\partition}{\pi}
\newcommand{\clu}{\mathsf{C}}
\newcommand{\cnum}{\mathsf{r}}
\newcommand{\pvec}{p}
\newcommand{\pmat}{\mathsf{P}}
\newcommand{\PTP}{\ensuremath{\myparen*{\pmat\tran \pmat}}}

\newcommand{\Ta}{T}
\newcommand{\Sa}{S}
\newcommand{\Tg}{\mathsf{T}}

\newcommand{\tran}{^{\operatorname{T}}}

\newcommand{\Ltwo}{\ensuremath{\cL_{2}}}
\newcommand{\Htwo}{\ensuremath{\cH_{2}}}
\newcommand{\Hinf}{\ensuremath{\cH_{\infty}}}

\DeclareMathOperator{\dif}{d\!}

\DeclarePairedDelimiter{\myparen}{\lparen}{\rparen}
\DeclarePairedDelimiter{\card}{\lvert}{\rvert}

\DeclarePairedDelimiterXPP{\normtwo}[1]{}{\lVert}{\rVert}{_{2}}{#1}
\DeclarePairedDelimiterXPP{\normF}[1]{}{\lVert}{\rVert}{_{\operatorname{F}}}{#1}
\DeclarePairedDelimiterXPP{\normLtwo}[1]{}{\lVert}{\rVert}{_{\Ltwo}}{#1}
\DeclarePairedDelimiterXPP{\normHtwo}[1]{}{\lVert}{\rVert}{_{\Htwo}}{#1}
\DeclarePairedDelimiterXPP{\normHinf}[1]{}{\lVert}{\rVert}{_{\Hinf}}{#1}
\DeclarePairedDelimiterXPP{\diag}[1]{\operatorname{diag}}{\lparen}{\rparen}{}{#1}
\DeclarePairedDelimiterXPP{\kernel}[1]{\operatorname{ker}}{\lparen}{\rparen}{}{#1}
\DeclarePairedDelimiterXPP{\image}[1]{\operatorname{im}}{\lparen}{\rparen}{}{#1}
\DeclarePairedDelimiterXPP{\col}[1]{\operatorname{col}}{\lparen}{\rparen}{}{#1}
\DeclarePairedDelimiterXPP{\eig}[1]{\sigma}{\lparen}{\rparen}{}{#1}
\DeclarePairedDelimiterXPP{\bigO}[1]{\mathcal{O}}{\lparen}{\rparen}{}{#1}

\let\le\leqslant%
\let\hat\widehat%
\let\bar\overline%

\title{Clustering-Based Model Order Reduction for Nonlinear Network Systems}
\author{%
  Peter Benner\thanks{%
    Max Planck Institute for Dynamics of Complex Technical Systems,
    Sandtorstr. 1,
    39106 Magdeburg,
    Germany,
    \texttt{\{benner,grundel,mlinaric\}@mpi-magdeburg.mpg.de}%
  } \and
  Sara Grundel\footnotemark[1] \and
  Petar Mlinari\'c\footnotemark[1]
}

\begin{document}
\maketitle

\begin{abstract}
  Clustering by projection has been proposed as a way to preserve network
  structure in linear multi-agent systems.
  Here, we extend this approach to a class of nonlinear network systems.
  Additionally, we generalize our clustering method which restores the network
  structure in an arbitrary reduced-order model obtained by projection.
  We demonstrate this method on a number of examples.
\end{abstract}

\section{Introduction}\label{sec:BGM20:intro}
Nonlinear network systems appear in various application areas, including energy
distribution networks, water networks, multi-robot networks, and chemical
reaction networks.
Model order reduction (MOR) enables faster simulation, optimization, and control
of large-scale network systems.
However, standard methods generally do not preserve the network structure.
Preserving the network structure is necessary, e.g., if an optimization method
assumes this structure.

Clustering was proposed in the literature as a way to preserve the multi-agent
structure.
Methods based on equitable partitions were described
in~\cite{MarEB08,RahJMetal09,morChaM11}
with an extension to almost equitable partitions in~\cite{MarEB10}.
Based on this, a priori error expressions were developed in~\cite{morMonTC14}
with generalizations in~\cite{morJonMGetal18}.
Ishizaki et al.~\cite{morIshKIetal14} developed a clustering-based $\Hinf$-MOR
method based on positive tridiagonalization and reducible clusters,
applicable to linear time-invariant systems with asymptotically stable and
symmetric dynamics matrices.
In~\cite{morIshKGetal15}, they presented an efficient clustering-based method
also based on reducible clusters for $\Htwo$-MOR of linear positive networks,
which include systems with Laplacian-based dynamics.
Cheng et al.~\cite{morCheKS16,morCheKS18} developed a method based on agent
dissimilarity.
Besselink et al.~\cite{morBesSJ16} studied networks of identical passive systems
over weighted and directed graphs with tree structures.

In this work, we extend the clustering-based approach for linear time-invariant
multi-agent systems from~\cite{morMliGB15,morMliGB16}.
There, we proposed a method combining
the iterative rational Krylov algorithm (IRKA)~\cite{morAntBG10} and
QR decomposition-based clustering~\cite{ZhaHDetal01}.
We generalize this approach to be able to combine
any projection-based MOR method and
clustering algorithm.
Extending to arbitrary projection-based MOR methods allows applying the method
to nonlinear network systems.
For the clustering algorithm,
we motivate the use of the k-means algorithm~\cite{HarW79}.
We show that for a class of nonlinear multi-agent systems,
clustering by Galerkin projection preserves network structure,
which additionally avoids the need for hyper-reduction to simplify the nonlinear
part.

The outline of this paper is as follows.
First, we provide some background information on linear multi-agent systems in
Section~\ref{sec:BGM20:prelim}.
In Section~\ref{sec:BGM20:clustering},
we recall our clustering-based MOR method for linear multi-agent systems and
generalize it to a framework which allows combining any projection-based MOR
method and clustering algorithm.
In Section~\ref{sec:BGM20:nonlinear-clustering},
we extend clustering by projection to a class of nonlinear multi-agent systems,
which also permits the applicability of our framework.
We demonstrate the approach numerically in Section~\ref{sec:BGM20:example} and
conclude with Section~\ref{sec:BGM20:con}.

We use $\imag$ to denote the imaginary unit ($\imag^2 = -1$),
$\bbC_-$ as the open left complex half-planes,
and $\bbC_+$ as the right.
Furthermore,
we use $\diag{v}$ to denote the diagonal matrix with the vector $v$ as its
diagonal and
$\col{v_1, v_2, \ldots, v_k}$ as the vector obtained by concatenating
$v_1, v_2, \ldots, v_k$.
We call a square matrix $A$ Hurwitz if all its eigenvalues have negative real
parts.
Similarly,
for square matrices $A$ and $B$,
with $B$ invertible,
we call the matrix pair $(A, B)$ Hurwitz if $B^{-1} A$ is Hurwitz.
For a rectangular matrix $A$,
$\image{A}$ denotes the subspace generated by the columns of $A$.
For a rational matrix function $\fundef{H}{\bbC}{\bbC^{p \times m}}$, i.e.,
a matrix-valued function whose components are rational functions,
the $\Htwo$ and $\Hinf$ norms are
\begin{align*}
  \normHtwo{H}
  & =
    \myparen*{
      \int_{-\infty}^\infty \normF{H(\imag \omega)}^2 \dif{\omega}
    }^{1/2}, \\
  \normHinf{H}
  & =
    \sup_{\omega \in \bbR} \normtwo{H(\imag \omega)},
\end{align*}
if all the poles of $H$ have negative real parts and undefined otherwise.

\section{Preliminaries}\label{sec:BGM20:prelim}
We present some basic concepts from graph theory in
Section~\ref{sec:BGM20:prelim-graph-theory},
graph partitions in Section~\ref{sec:BGM20:prelim-partitions},
before moving on to linear multi-agent systems in
Section~\ref{sec:BGM20:prelim-mas} and
clustering-based MOR in Section~\ref{sec:BGM20:prelim-clust-mor}.
Additionally, we give remarks on MOR for non-asymptotically stable linear
multi-agent systems in Section~\ref{sec:BGM20:prelim-unstable-mor}.

\subsection{Graph theory}\label{sec:BGM20:prelim-graph-theory}
The notation in this section is based on~\cite{MesE10} and~\cite{GodR01}.

A \emph{graph} $\graph$ consists of a \emph{vertex set} $\vset$ and an
\emph{edge set} $\eset$ encoding the relation between vertices.
\emph{Undirected} graphs are those for which the edge set is a subset of the
set of all unordered pairs of vertices, i.e.,
$\eset \subseteq \{\{i, j\} : i, j \in \vset,\ i \neq j\}$.
On the other hand, a graph is \emph{directed} if
$\eset \subseteq \{(i, j) : i, j \in \vset,\ i \neq j\}$.
We think of an edge~$(i, j)$ as an arrow starting from vertex~$i$ and ending
at~$j$.
We only consider \emph{simple} graphs, i.e.,
graphs without self-loops or multiple copies of the same edge.
Additionally, we only consider \emph{finite graphs}, i.e.,
graphs with a finite number of vertices $\vnum \coloneqq \card{\vset}$.
Without loss of generality,
let $\vset = \{1, 2, \dots, \vnum\}$.

For an undirected graph, a \emph{path} of length $\ell$ is a sequence of
distinct vertices $i_0, i_1, \dots, i_{\ell}$ such that
$\{i_k, i_{k + 1}\} \in \eset$ for $k = 0, 1, \dots, \ell - 1$.
For a directed graph, a \emph{directed path} of length $\ell$ is a sequence of
distinct vertices $i_0, i_1, \dots, i_{\ell}$ such that
$(i_k, i_{k + 1}) \in \eset$ for $k = 0, 1, \dots, \ell - 1$.
An undirected graph is \emph{connected} if there is a path between any two
distinct vertices $i, j \in \vset$.
A directed graph is \emph{strongly connected} if there is a directed path
between any two distinct vertices $i, j \in \vset$.

We can associate weights to edges of a graph by a \emph{weight function}
$\fundef{\wfun}{\eset}{\bbR}$.
If $\wfun(\edge) > 0$ for all $\edge \in \eset$, the tuple
$\graph = (\vset, \eset, \wfun)$ is called a \emph{weighted graph}.
In the following, we will focus on weighted graphs.
In particular, we will directly generalize concepts for unweighted graphs
from~\cite{MesE10,GodR01}, as was done in~\cite{morMonTC14}.

The \emph{adjacency matrix}
$\amat = {[\weight_{ij}]}_{i, j \in \vset} \in \bbR^{\vnum \times \vnum}$ of an
undirected weighted graph is defined component-wise by
\begin{align*}
  \weight_{ij} \coloneqq
  \begin{cases}
    \wfun(\{i, j\}), & \text{if } \{i, j\} \in \eset, \\
    0, & \text{otherwise},
  \end{cases}
\end{align*}
and for a directed weighted graph as
\begin{align*}
  \weight_{ij} \coloneqq
  \begin{cases}
    \wfun((j, i)), & \text{if } (j, i) \in \eset, \\
    0, & \text{otherwise}.
  \end{cases}
\end{align*}
For every vertex $i \in \vset$, its \emph{in-degree} is
$\degree_i \coloneqq \sum_{j = 1}^{\vnum}{\weight_{ij}}$.
The diagonal matrix
$\dmat \coloneqq \diag{\degree_1, \degree_2, \ldots, \degree_{\vnum}}$ is called
the \emph{in-degree matrix}.
Notice that $\dmat = \diag{\amat \one}$, where $\one$ is the vector of all ones.

Let $\edge_1, \edge_2, \dots, \edge_{\card{\eset}}$ be all the edges of
$\graph$ in some order.
The \emph{incidence matrix} $\rmat \in \bbR^{\vnum \times \card{\eset}}$ of a
directed graph $\graph$ is defined component-wise
\begin{align*}
  {[\rmat]}_{ik} \coloneqq
  \begin{cases}
    -1, & \text{if } \edge_k = (i, j) \text{ for some } j \in \vset, \\
    1, & \text{if } \edge_k = (j, i) \text{ for some } j \in \vset, \\
    0, & \text{otherwise}.
  \end{cases}
\end{align*}
If $\graph$ is undirected, we assign some orientation to every edge to define a
directed graph $\graph^o$, and define the incidence matrix of $\graph$ to be the
incidence matrix of $\graph^o$.
The \emph{weight matrix} is defined as
$\wmat \coloneqq
\diag{\wfun(\edge_1), \wfun(\edge_2), \dots, \wfun(\edge_{\card{\eset}})}$.

The \emph{(in-degree) Laplacian matrix} $\lmat$ is defined by
$\lmat \coloneqq \dmat - \amat$.
For undirected graphs, it can be checked that $\lmat = \rmat \wmat \rmat\tran$,
using
\begin{align*}
  \rmat \wmat \rmat\tran
  = \sum_{\{i, j\} \in \eset}{\weight_{ij} (e_i - e_j) (e_i - e_j)\tran},
\end{align*}
which is independent of the order of edges defining $\rmat$ and $\wmat$ or the
orientation of edges in $\graph^o$.
From the definition of $\lmat$, it directly follows that the sum of each row in
$\lmat$ is zero, i.e., $\lmat \one = 0$.
From $\lmat = \rmat \wmat \rmat\tran$, we immediately see that, for undirected
weighted graphs, the Laplacian matrix $\lmat$ is symmetric positive
semidefinite.

The following theorem,
based on Theorem~2.8 in~\cite{MesE10},
states how connectedness of a graph is related to the spectral properties of
$\lmat$.
\begin{theorem}
  Let $\graph = (\vset, \eset, \wfun)$ be an undirected weighted graph, $\lmat$
  its Laplacian matrix, and
  $0 = \lambda_1 \le \lambda_2 \le \dots \le \lambda_{\vnum}$ the eigenvalues of
  $\lmat$.
  Then the following statements are equivalent:
  \begin{enumerate}
  \item $\graph$ is connected,
  \item $\lambda_2 > 0$,
  \item $\kernel{\lmat} = \image{\one}$.
  \end{enumerate}
\end{theorem}

\subsection{Graph partitions}\label{sec:BGM20:prelim-partitions}
A nonempty subset $\clu \subseteq \vset$ is called a \emph{cluster} of $\vset$.
A \emph{graph partition} $\pi$ of the graph $\graph$ is a partition of its
vertex set $\vset$.
The \emph{characteristic vector} of a cluster $\clu \subseteq \vset$ is the
vector $\pvec(\clu) \in \bbR^{\vnum}$ defined with
\begin{equation*}
  {[\pvec(\clu)]}_i \coloneqq
  \begin{cases}
    1 & \text{if } i \in \clu, \\
    0 & \text{otherwise}.
  \end{cases}
\end{equation*}
The \emph{characteristic matrix of a partition}
$\partition = \{\clu_1, \clu_2, \ldots, \clu_\cnum\}$ is the matrix
$\pmat \in \bbR^{\vnum \times \cnum}$ defined by
\begin{equation*}
  \pmat \coloneqq
  \begin{bmatrix}
    \pvec(\clu_1) & \pvec(\clu_2) & \cdots & \pvec(\clu_{\cnum})
  \end{bmatrix}.
\end{equation*}
Note that
$\pmat\tran \pmat
= \diag{\card{\clu_1}, \card{\clu_2}, \ldots, \card{\clu_\cnum}}$.

\subsection{Linear multi-agent systems}\label{sec:BGM20:prelim-mas}
Here, we focus on linear time-invariant multi-agent systems
(cf.~\cite{morBesSJ16,morCheKS16,morCheKS18,morIshKIetal14,morIshKGetal15,morIshKI16,morMonTC13,morMonTC14}).
Additionally, we restrict ourselves to multi-agent systems
defined over an undirected, weighted, and connected graph
$\graph = (\vset, \eset, \wfun)$.

The dynamics of the $i$th agent, for
$i \in \vset = \{1, 2, \dots, \vnum\}$, is
\begin{align*}
  E \dot{x}_i(t) & = A x_i(t) + B v_i(t), \\
  z_i(t) & = C x_i(t),
\end{align*}
with system matrices $E, A \in \Rnn$, input matrix $B \in \Rnm$, output matrix
$C \in \Rpn$, state $x_i(t) \in \Rn$, input $v_i(t) \in \Rm$,
and output $z_i(t) \in \Rp$.
We assume the matrix $E$ to be invertible.
The interconnections are
\begin{align*}
  \mass_i v_i(t)
  & =
    \sum_{j = 1}^{\vnum}{\weight_{ij} K \myparen*{z_j(t) - z_i(t)}}
    + \sum_{k = 1}^{\inum}{\iw_{ik} u_k(t)},
\end{align*}
for $i = 1, 2, \dots, \vnum$,
with inertias $\mass_i > 0$,
coupling matrix $K \in \Rmp$,
external inputs $u_k(t) \in \Rm$,
$k = 1, 2, \dots, \inum$,
where $\amat = [\weight_{ij}]$ is the adjacency matrix of the graph $\graph$.
The outputs are
\begin{align*}
  y_{\ell}(t)
  & =
    \sum_{j = 1}^{\vnum}{\ow_{\ell j} z_j(t)}
\end{align*}
for $\ell = 1, 2, \dots, \onum$.
Define
\begin{gather*}
  \mmat \coloneqq \diag{\mass_i} \in \bbR^{\vnum \times \vnum},\
  \imat \coloneqq {[\iw_{ik}]} \in \bbR^{\vnum \times \inum},\
  \omat \coloneqq {[\ow_{\ell j}]} \in \bbR^{\onum \times \vnum}, \\
  x(t) \coloneqq \col{x_i(t)} \in \bbR^{\vnum n},\
  v(t) \coloneqq \col{v_i(t)} \in \bbR^{\vnum m},\
  z(t) \coloneqq \col{z_i(t)} \in \bbR^{\vnum p}, \\
  u(t) \coloneqq \col{u_k(t)} \in \bbR^{\inum m},
  \text{ and }
  y(t) \coloneqq \col{y_{\ell}(t)} \in \bbR^{\onum p}.
\end{gather*}
Then the agent dynamics can be rewritten as
\begin{align*}
  (I_{\vnum} \otimes E) \dot{x}(t)
  & =
    (I_{\vnum} \otimes A) x(t)
    + (I_{\vnum} \otimes B) v(t), \\
  z(t)
  & =
    (I_{\vnum} \otimes C) x(t),
\end{align*}
interconnection as
\begin{align*}
  (\mmat \otimes I_n) v(t)
  & =
    (-\lmat \otimes K) z(t)
    + (\imat \otimes I_m) u(t),
\end{align*}
and output as
\begin{align*}
  y(t)
  & =
    (\omat \otimes I_p) z(t).
\end{align*}
Therefore, we have
\begin{equation}\label{eq:BGM20:general-mas}
  \begin{aligned}
    (\mmat \otimes E) \dot{x}(t)
    & =
      (\mmat \otimes A - \lmat \otimes B K C) x(t)
      + (\imat \otimes B) u(t), \\
    y(t)
    & =
      (\omat \otimes C) x(t).
  \end{aligned}
\end{equation}

Of particular interest are \emph{leader-follower multi-agent systems} where only
some agents (\emph{leaders}) receive external input, while other agents
(\emph{followers}) receive no inputs.
Let $\inum \in \{1, 2, \dots, \vnum\}$ be the number of leaders,
$\VL = \{\ver_1, \ver_2, \dots, \ver_{\inum}\} \subseteq \vset$ the set of
leaders, and $\VF = \vset \setminus \VL$ the set of followers.
Then, with $\imat$ defined by
\begin{align*}
	\iw_{ik} \coloneqq
	\begin{cases}
		1, & \text{if } i = \ver_k, \\
		0, & \text{otherwise},
	\end{cases}
\end{align*}
the system~\eqref{eq:BGM20:general-mas} becomes a leader-follower multi-agent
system.
One important class are multi-agent systems with \emph{single-integrator
  agents}, i.e., with $n = 1$, $A = 0$, and $B = C = K = E = 1$.
Thus, system~\eqref{eq:BGM20:general-mas} becomes
\begin{equation}\label{eq:BGM20:si_mas}
  \begin{aligned}
    \mmat \dot{x}(t) & = -\lmat x(t) + \imat u(t), \\
    y(t) & = \omat x(t).
  \end{aligned}
\end{equation}

The property of interest for multi-agent systems is \emph{synchronization}.
\begin{definition}
  The system~\eqref{eq:BGM20:general-mas} is \emph{synchronized} if
  \begin{align*}
    \lim_{t \to \infty}{(x_i(t) - x_j(t))} = 0,
  \end{align*}
  for all $i, j \in \vset$ and all initial conditions $x(0) = x_0$
  and $u \equiv 0$.
\end{definition}
In words, this means that the agents' states converge to the same trajectory for
zero input and arbitrary initial condition.
The following results gives a characterization
(\cite[Theorem~1]{LiDCetal10},~\cite[Lemma~4.2]{morMonTC13}).
\begin{proposition}\label{thm:BGM20:mas-sync-equiv}
  Let a system~\eqref{eq:BGM20:general-mas} be given,
  where $\lmat$ is the Laplacian matrix of an undirected, weighted, and
  connected graph.
  Then the system~\eqref{eq:BGM20:general-mas} is synchronized if and only if
  $(A - \lambda B K C, E)$ is Hurwitz
  for all nonzero eigenvalues $\lambda$ of $(\lmat, \mmat)$.
\end{proposition}
Note that linear multi-agent systems with single integrator agents,
as in~\eqref{eq:BGM20:si_mas},
are always synchronized since $(A - \lambda B K C, E) = (-\lambda, 1)$.

\subsection{Clustering-based model order
  reduction}\label{sec:BGM20:prelim-clust-mor}
By choosing some matrices $V, W \in \bbR^{\vnum n \times \cnum n}$, we get the
reduced model for~\eqref{eq:BGM20:si_mas}
\begin{equation}\label{eq:BGM20:red_pg}
	\begin{aligned}
		W\tran \mmat V \dot{\hx}(t)
    &
      = -W\tran \lmat V \hx(t) + W\tran \imat u(t), \\
		\hy(t)
    &
      = \omat V \hx(t),
	\end{aligned}
\end{equation}
or, for~\eqref{eq:BGM20:general-mas},
\begin{equation}\label{eq:BGM20:general-reg-pg}
  \begin{aligned}
    W\tran (\mmat \otimes E) V \dot{\hx}(t)
    & =
      W\tran (\mmat \otimes A - \lmat \otimes B K C) V \hx(t)
      + W\tran (\imat \otimes B) u(t), \\
    \hy(t)
    & =
      (\omat \otimes C) V \hx(t),
  \end{aligned}
\end{equation}
which is not necessarily a multi-agent system.
As suggested in~\cite{morCheKS16} (similar to~\cite{morMonTC14,morIshKIetal14}),
using
\begin{equation}\label{eq:BGM20:VWP}
  V
  = W
  = \pmat,
\end{equation}
in~\eqref{eq:BGM20:red_pg},
or in general
\begin{equation}\label{eq:BGM20:VWPn}
  V
  = W
  = \pmat \otimes I_n,
\end{equation}
in~\eqref{eq:BGM20:general-reg-pg},
preserves the structure,
where $\pmat$ is a characteristic matrix of a partition $\partition$ of the
vertex set $\vset$.
In particular,
$\pmat\tran \mmat \pmat$ is a positive definite diagonal matrix and
$\pmat\tran \lmat \pmat$ is the Laplacian matrix of the reduced graph.

\subsection{Model reduction for non-asymptotically stable
  systems}\label{sec:BGM20:prelim-unstable-mor}
Note that the system~\eqref{eq:BGM20:si_mas} is not (internally) asymptotically
stable since $\lmat$ has a zero eigenvalue.
Similarly, the system~\eqref{eq:BGM20:general-mas} is not asymptotically stable
if $A$ is not Hurwitz.
First, we discuss a decomposition into the asymptotically and the
non-asymptotically stable part.
This allows an extension of MOR methods and the computation of system norms.
Next, we analyze stability of clustering-based reduced models.

For an arbitrary linear time-invariant system
\begin{align*}
  \cE \dot{x}(t) & = \cA x(t) + \cB u(t), \\
  y(t) & = C x(t),
\end{align*}
with invertible $\cE$,
let
$\cT =
\begin{bmatrix}
  \cT_- & \cT_+
\end{bmatrix}
$
and
$\cS =
\begin{bmatrix}
  \cS_- & \cS_+
\end{bmatrix}
$
be invertible matrices such that
\begin{align*}
  \cS\tran \cE \cT =
  \begin{bmatrix}
    \cE_- & 0 \\
    0 & \cE_+
  \end{bmatrix}, \quad
  \cS\tran \cA \cT =
  \begin{bmatrix}
    \cA_- & 0 \\
    0 & \cA_+
  \end{bmatrix},
\end{align*}
with
$\eig{\cA_-, \cE_-} \subset \bbC_-$ and
$\eig{\cA_+, \cE_+} \subset \bar{\bbC_+}$.
In particular,
this means that $\image{\cT_-}$ is a direct sum of generalized (right)
eigenspaces corresponding to the eigenvalues of $(\cA, \cE)$ with negative real
parts and analogously for $\image{\cT_+}, \image{\cS_-}, \image{\cS_+}$.
If we denote
\begin{align*}
  \cS\tran \cB =
  \begin{bmatrix}
    \cB_- \\
    \cB_+
  \end{bmatrix}, \quad
  \cC \cT =
  \begin{bmatrix}
    \cC_- & \cC_+
  \end{bmatrix},
\end{align*}
this gives us that $\cH = \cH_- + \cH_+$,
where
\begin{align*}
  \cH(s)
  & =
    \cC \myparen*{s \cE - \cA}^{-1} \cB,
    \displaybreak[0] \\
  \cH_-(s)
  & =
    \cC_- \myparen*{s \cE_- - \cA_-}^{-1} \cB_-,
    \displaybreak[0] \\
  \cH_+(s)
  & =
    \cC_+ \myparen*{s \cE_+ - \cA_+}^{-1} \cB_+.
\end{align*}
Note that $\cH_-$ and $\cH_+$ have poles in $\bbC_-$ and $\bar{\bbC_+}$,
respectively.
For the transfer function $\hcH$ of a reduced model to be such that
$\normHtwo{\cH - \hcH}$ and $\normHinf{\cH - \hcH}$ are defined,
it is necessary that $\cH$ and $\hcH$ have the same non-asymptotically stable
part.
This means that $\hcH = \hcH_- + \cH_+$ for some $\hcH_-$,
$\hcH_-(s) = \hcC_- \myparen{s \hcE_- - \hcA_-}^{-1} \hcB_-$,
with poles in $\bbC_-$, i.e.,
$\hcH_-$ is a reduced model for $\cH_-$.
If we use a projection-based MOR method with matrices $\cV_-, \cW_-$ to get
\begin{align*}
  \hcE_- = \cW_-\tran \cE_- \cV_-, \quad
  \hcA_- = \cW_-\tran \cA_- \cV_-, \quad
  \hcB_- = \cW_-\tran \cB_-, \quad
  \hcC_- = \cC_- \cV_-,
\end{align*}
then the overall basis matrices are
\begin{equation}\label{eq:BGM20:pg_stab}
  \cV =
  \begin{bmatrix}
    \cT_- \cV_- & \cT_+
  \end{bmatrix}
  \quad \text{and} \quad
  \cW =
  \begin{bmatrix}
    \cS_- \cW_- & \cS_+
  \end{bmatrix}.
\end{equation}
Then we can compute the norm of $\cH - \hcH$ by computing the norm of
$\cH_- - \hcH_-$.

To analyze linear multi-agent systems,
we want to find an invertible matrix $\Tg$ such that
\begin{align*}
  \Tg\tran \mmat \Tg =
  \begin{bmatrix}
    \mmat_- & 0 \\
    0 & \mass_+
  \end{bmatrix}
  \text{ and }
  \Tg\tran \lmat \Tg =
  \begin{bmatrix}
    \lmat_- & 0 \\
    0 & 0
  \end{bmatrix},
\end{align*}
where $\eig{-\lmat_-, \mmat_-} \subset \bbC_-$.
We see that if
\begin{align*}
  \Tg =
  \begin{bmatrix}
    \Tg_- & \one_{\vnum} \\
  \end{bmatrix},
\end{align*}
then
\begin{align*}
  \Tg\tran \mmat \Tg =
  \begin{bmatrix}
    \Tg_-\tran \mmat \Tg_-
    & \Tg_-\tran \mmat \one_{\vnum} \\
    \one_{\vnum}\tran \mmat \Tg_-
    & \one_{\vnum}\tran \mmat \one_{\vnum}
  \end{bmatrix}
  \text{ and }
  \Tg\tran \lmat \Tg =
  \begin{bmatrix}
    \Tg_-\tran \lmat \Tg_- & 0 \\
    0 & 0
  \end{bmatrix}.
\end{align*}
To have $\Tg_-\tran \mmat \one_{\vnum} = 0$, we need for the columns of $\Tg_-$
to be orthogonal to $\mmat \one_{\vnum}$.
This will also ensure that
$
\eig{-\lmat_-, \mmat_-}
= \eig{-\Tg_-\tran \lmat \Tg_-, \Tg_-\tran \mmat \Tg_-}
\subset \bbC_-
$.
Additionally, $\Tg_-$ should be such that both $\Tg_-\tran \mmat \Tg_-$ and
$\Tg_-\tran \lmat \Tg_-$ are sparse.
We chose the form
\begin{align*}
  \Tg_- =
  \begin{bmatrix}
    \alpha_1 \\
    -\beta_1 & \alpha_2 \\
    & -\beta_2 & \ddots \\
    & & \ddots & \alpha_{\vnum - 1} \\
    & & & -\beta_{\vnum - 1}
  \end{bmatrix}
\end{align*}
with some $\alpha_i, \beta_i > 0$, $i = 1, 2, \dots, \vnum - 1$,
which we determine next.
From $e_i\tran \Tg_-\tran \mmat \one_{\vnum} = 0$, we find
$\alpha_i \mass_i = \beta_i \mass_{i + 1}$.
If we additionally set $\alpha_i^2 + \beta_i^2 = 1$, we get
\begin{align*}
  \alpha_i = \frac{\mass_{i + 1}}{\sqrt{\mass_i^2 + \mass_{i + 1}^2}}
  \quad \text{and} \quad
  \beta_i = \frac{\mass_i}{\sqrt{\mass_i^2 + \mass_{i + 1}^2}}.
\end{align*}
Therefore, the decomposition of a multi-agent system~\eqref{eq:BGM20:si_mas} is
\begin{align*}
  H(s)
  & =
    \omat \Tg_- \myparen*{s \mmat_- - \lmat_-}^{-1} \Tg_-\tran \imat
    + \frac{1}{s \mass_+} \omat \one_{\vnum} \one_{\vnum}\tran \imat
\end{align*}
Similarly,
for the reduced model~\eqref{eq:BGM20:red_pg} with~\eqref{eq:BGM20:VWP},
we get that the non-asymptotically stable part is
\begin{align*}
  \frac{1}{s \one_{\cnum}\tran \pmat\tran \mmat \pmat \one_{\cnum}}
  \omat \pmat \one_{\cnum} \one_{\cnum}\tran \pmat\tran \imat,
\end{align*}
which is equal to the non-asymptotically stable part of the original model since
$\pmat \one_{\cnum} = \one_{\vnum}$.
Therefore, the transfer function of the error system has only poles with
negative real parts.

Next,
to analyze system~\eqref{eq:BGM20:general-mas},
let $\Ta$ and $\Sa$ be invertible matrices such that
\begin{align*}
  \Sa\tran E \Ta =
  \begin{bmatrix}
    E_- & 0 \\
    0 & E_+
  \end{bmatrix}
  \text{ and }
  \Sa\tran A \Ta =
  \begin{bmatrix}
    A_- & 0 \\
    0 & A_+
  \end{bmatrix},
\end{align*}
where $\eig{A_-, E_-} \subset \bbC_-$ and
$\eig{A_+, E_+} \subset \bar{\bbC_+}$.
Then
\begin{align*}
  &
    \begin{bmatrix}
      I_{(\vnum - 1) n} & 0 \\
      0 & \Sa
    \end{bmatrix}\tran
    (\Tg \otimes I_n)\tran
    (\mmat \otimes E)
    (\Tg \otimes I_n)
    \begin{bmatrix}
      I_{(\vnum - 1) n} & 0 \\
      0 & \Ta
    \end{bmatrix}
    \displaybreak[0]  \\
  & =
    \begin{bmatrix}
      I_{(\vnum - 1) n} & 0 \\
      0 & \Sa
    \end{bmatrix}\tran
    \begin{bmatrix}
      \mmat_- \otimes E & 0 \\
      0 & \mass_+ E
    \end{bmatrix}
    \begin{bmatrix}
      I_{(\vnum - 1) n} & 0 \\
      0 & \Ta
    \end{bmatrix}
    \displaybreak[0] \\
  & =
    \begin{bmatrix}
      \mmat_- \otimes E & 0 & 0 \\
      0 & \mass_+ E_- & 0 \\
      0 & 0 & \mass_+ E_+
    \end{bmatrix}
\end{align*}
and
\begin{align*}
  &
    \begin{bmatrix}
      I_{(\vnum - 1) n} & 0 \\
      0 & \Sa
    \end{bmatrix}\tran
    (\Tg \otimes I_n)\tran
    (\mmat \otimes A - \lmat \otimes B K C)
    (\Tg \otimes I_n)
    \begin{bmatrix}
      I_{(\vnum - 1) n} & 0 \\
      0 & \Ta
    \end{bmatrix}
    \displaybreak[0] \\
  & =
    \begin{bmatrix}
      I_{(\vnum - 1) n} & 0 \\
      0 & \Sa
    \end{bmatrix}\tran
    \begin{bmatrix}
      \mmat_- \otimes A - \lmat_- \otimes B K C & 0 \\
      0 & \mass_+ A
    \end{bmatrix}
    \begin{bmatrix}
      I_{(\vnum - 1) n} & 0 \\
      0 & \Ta
    \end{bmatrix}
    \displaybreak[0] \\
  & =
    \begin{bmatrix}
      \mmat_- \otimes A - \lmat_- \otimes B K C & 0 & 0 \\
      0 & \mass_+ A_- & 0 \\
      0 & 0 & \mass_+ A_+
    \end{bmatrix}.
\end{align*}
Since the original system is assumed to be synchronized,
we have that
$\eig{\mmat_- \otimes A - \lmat_- \otimes B K C, \mmat_- \otimes E}
\subset \bbC_-$.
Note that the overall transformation matrices are
\begin{align*}
  \cT
  & =
    (\Tg \otimes I_n)
    \begin{bmatrix}
      I_{(\vnum - 1) n} & 0 \\
      0 & \Ta
    \end{bmatrix}
    =
    \myparen*{
      \begin{bmatrix}
        \Tg_- & \one_{\vnum}
      \end{bmatrix}
      \otimes I_n
    }
    \begin{bmatrix}
      I_{(\vnum - 1) n} & 0 \\
      0 & \Ta
    \end{bmatrix}
    \displaybreak[0] \\
  & =
    \begin{bmatrix}
      \Tg_- \otimes I_n & \one_{\vnum} \otimes I_n
    \end{bmatrix}
    \begin{bmatrix}
      I_{(\vnum - 1) n} & 0 \\
      0 & \Ta
    \end{bmatrix}
    =
    \begin{bmatrix}
      \Tg_- \otimes I_n & (\one_{\vnum} \otimes I_n) \Ta
    \end{bmatrix}
    \displaybreak[0] \\
  & =
    \begin{bmatrix}
      \Tg_- \otimes I_n & \one_{\vnum} \otimes \Ta
    \end{bmatrix}
    =
    \begin{bmatrix}
      \Tg_- \otimes I_n & \one_{\vnum} \otimes \Ta_- & \one_{\vnum} \otimes \Ta_+
    \end{bmatrix},
    \displaybreak[0] \\
  \cS
  & =
    \begin{bmatrix}
      \Tg_- \otimes I_n & \one_{\vnum} \otimes \Sa_- & \one_{\vnum} \otimes \Sa_+
    \end{bmatrix}.
\end{align*}
Therefore, the non-asymptotically stable part is
\begin{align*}
  &
    (\omat \otimes C) (\one_{\vnum} \otimes \Ta_+)
    \myparen*{s \mass_+ E_+ - \mass_+ A_+}^{-1}
    (\one_{\vnum} \otimes \Sa_+)\tran (\imat \otimes B) \\
  & =
    \myparen*{\omat \one_{\vnum} \otimes C \Ta_+}
    \myparen*{s \mass_+ E_+ - \mass_+ A_+}^{-1}
    \myparen*{\one_{\vnum}\tran \imat \otimes \Sa_+\tran B}
\end{align*}
Assuming that a clustering-based reduced model is synchronized,
we see that it has the same non-asymptotically stable part as the original
model.

It remains to consider synchronization preservation.
In the single-integrator case, clustering using any partition preserves
synchronization.
In the general case, using Theorem~\ref{thm:BGM20:mas-sync-equiv}, we need that
$(A - \hat{\lambda} B K C, E)$ is Hurwitz
for all nonzero eigenvalues $\hat{\lambda}$ of $(\hat{\lmat}, \hat{\mmat})$.
Since in general $\eig{\hat{\lmat}, \hat{\mmat}}$ is not a subset of
$\eig{\lmat, \mmat}$,
we need an additional assumption.
Based on the interlacing property~\cite{GolV13},
we know that all nonzero eigenvalues $(\hat{\lmat}, \hat{\mmat})$ are in
$[\lambda_2, \lambda_{\vnum}]$,
where $0 = \lambda_1 < \lambda_2 \le \dots \le \lambda_{\vnum}$ are the
eigenvalues of $(\lmat, \mmat)$.
Therefore,
if $(A - \lambda B K C, E)$ is Hurwitz
for all $\lambda \in [\lambda_2, \lambda_{\vnum}]$,
we get that every partition preserves synchronization.

\section{Clustering for linear multi-agent systems}\label{sec:BGM20:clustering}
In this section,
we motivate our general approach for clustering-based linear multi-agent
systems.
Since clustering is generally a difficult combinatorial problem (see,
e.g.,~\cite{Sch07}),
we propose a heuristic approach for finding suboptimal partitions.

For simplicity,
we first consider multi-agent systems with single-integrator agents as
in~\eqref{eq:BGM20:si_mas}.
Let
\begin{align*}
	H(s)
  & =
    \omat {(s \mmat + \lmat)}^{-1} \imat, \\
	\hH(s)
  & =
    \omat V \myparen*{s W\tran \mmat V + W\tran \lmat V}^{-1} W\tran \imat
\end{align*}
be the transfer functions of systems~\eqref{eq:BGM20:si_mas}
and~\eqref{eq:BGM20:red_pg},
respectively,
where $V, W \in \bbR^{\vnum \times r_P}$ are obtained using a projection-based
method such as balanced truncation or IRKA and the construction
in~\eqref{eq:BGM20:pg_stab}.

In~\cite{morMliGB15},
motivated by~\eqref{eq:BGM20:VWP} and
the properties of clustering using QR decomposition with column pivoting
(see~\cite[Section~3]{ZhaHDetal01},~\cite[Lemma~1]{morMliGB15}),
we proposed applying it to the set of rows of $V$ or $W$ to recover the
partition.
Here, we want to emphasize that the approach is not restricted to this choice of
clustering algorithm.
In particular,
the following result on the forward error in the Petrov-Galerkin projection
(\cite[Theorem~3.3]{morBeaGW12})
motivates using the k-means clustering~\cite{HarW79}.
\begin{theorem}\label{thm:BGM:bound}
  Let $V_1, V_2, W_1, W_2 \in \bbR^{\vnum \times r_P}$ be full-rank matrices and
  \begin{align*}
    \cV_i = \image{V_i}, \quad
    \cW_i = \image{W_i}, \quad
    \hH_i(s)
    = C V_i \myparen*{s W_i\tran E V_i - W_i\tran A V_i}^{-1} W_i\tran B,
  \end{align*}
  for $i = 1, 2$.
  Then
  \begin{align*}
    \frac{\normHinf[\big]{\hH_1 - \hH_2}}{
      \frac{1}{2} \myparen*{\normHinf[\big]{\hH_1} + \normHinf[\big]{\hH_2}}}
    & \le
      M \max(\sin{\Theta(\cV_1, \cV_2)}, \sin{\Theta(\cW_1, \cW_2)}),
  \end{align*}
  where
  \begin{align*}
    M
    & =
      2 \max(M_1, M_2), \\
    M_1
    & =
      \frac{\max_{\omega \in \bbR}
        \normtwo{C}
        \normtwo*{V_1
          \myparen*{\imag \omega W_1\tran E V_1 - W_1\tran A V_1}^{-1}
          W_1\tran B}
        \normtwo*{\hH_1(\imag \omega)}^{-1}
      }{\min_{\omega \in \bbR} \cos
        \Theta\myparen*{
          \kernel*{W_2\tran \myparen{\imag \omega E - A}^{-1}}^{\perp},
          \cV_2
        }}, \\
    M_2
    & =
      \frac{\max_{\omega \in \bbR}
        \normtwo*{C V_2
          \myparen*{\imag \omega W_2\tran E V_2 - W_2\tran A V_2}^{-1}
          W_2\tran}
        \normtwo{B}
        \normtwo*{\hH_2(\imag \omega)}^{-1}
      }{\min_{\omega \in \bbR} \cos
        \Theta\myparen*{\image*{\myparen{\imag \omega E - A}^{-1} V_1}, \cW_1}},
  \end{align*}
  and $\Theta(\cM, \cN)$ is the largest principal angle between subspaces
  $\cM, \cN \subseteq \Rn$.
\end{theorem}
The motivation for looking at this bound is,
if we take $\hH_1$ to be a projection-based reduced-order model that is very
close to the original model, i.e., $\normHinf{H - \hH_1}$ is small,
then we could look for a clustering-based reduced-order model $\hH_2$ by finding
a characteristic matrix of a partition $\pmat$ such that $\image{\pmat}$ is
close to $\cV_1$ and $\cW_1$.

Note that,
to use Theorem~\ref{thm:BGM:bound},
we need to use the asymptotically stable parts of $\hH_1$ and $\hH_2$ such that
$\normHinf{\hH_1}$ and $\normHinf{\hH_2}$ are defined.
As discussed in Section~\ref{sec:BGM20:prelim-unstable-mor},
$\hH_1$ and $\hH_2$ need to have the same non-asymptotically stable part as $H$
for the $\Htwo$ and $\Hinf$ errors to be defined.

Next, we show how the bounds motivate the use of k-means clustering.
The sine of the largest principal angle between two subspaces
$\cV_1, \cV_2 \subseteq \Rn$ is defined by (see~\cite[Section~3.1]{morBeaGW12})
\begin{align*}
  \sin{\Theta(\cV_1, \cV_2)}
  \coloneqq \sup_{v_1 \in \cV_1} \inf_{v_2 \in \cV_2}
    \frac{\normtwo{v_2 - v_1}}{\normtwo{v_1}}.
\end{align*}
Furthermore,
if $\Pi_1$ and $\Pi_2$ are orthogonal projectors onto $\cV_1$ and $\cV_2$,
then $\sin{\Theta(\cV_1, \cV_2)} = \normtwo{(I - \Pi_2) \Pi_1}$.
Therefore, we have
\begin{align*}
  \sin{\Theta(\cV_1, \cV_2)}
  = \normtwo*{\myparen*{I - V_2 V_2\tran} V_1},
\end{align*}
for any $V_1$ and $V_2$ with orthonormal columns such that
$\cV_1 = \image{V_1}$, $\cV_2 = \image{V_2}$.
If additionally
$V_1$ is the $V \in \bbR^{\vnum \times r_P}$ from the projection-based method and
$V_2 = \pmat \PTP^{-1/2} \in \bbR^{\vnum \times \cnum}$,
then
\begin{align*}
  \myparen*{\sin{\Theta(\cV_1, \cV_2)}}^2
  & \le
    \normF*{\myparen*{I - \pmat \PTP^{-1} \pmat\tran} V}^2
    \displaybreak[0] \\
  & =
    \normF*{\myparen*{I -
      \begin{bmatrix}
        \pvec(\clu_1) & \cdots & \pvec(\clu_{\cnum})
      \end{bmatrix}
      \begin{bmatrix}
        \card{\clu_1}^{-1} \\
        & \ddots \\
        & & \card{\clu_{\cnum}}^{-1}
      \end{bmatrix}
      \begin{bmatrix}
        \pvec(\clu_1)\tran \\
        \vdots \\
        \pvec(\clu_{\cnum})\tran
      \end{bmatrix}}
    V}^2
    \displaybreak[0] \\
  & =
    \normF*{\myparen*{I - \sum_{i = 1}^{\cnum}
      \frac{1}{\card{\clu_i}} \pvec(\clu_i) \pvec(\clu_i)\tran} V}^2
    \displaybreak[0] \\
  & =
    \normF*{
      \sum_{i = 1}^{\cnum} \sum_{p \in \clu_i}
      \myparen*{e_p e_p\tran
        - \frac{1}{\card{\clu_i}} e_p \pvec(\clu_i)\tran}
      V
    }^2
    \displaybreak[0] \\
  & =
    \sum_{i = 1}^{\cnum} \sum_{p \in \clu_i}
    \normtwo*{V_{p, :} - \frac{1}{\card{\clu_i}} \sum_{q \in \clu_i} V_{q, :}}^2,
\end{align*}
which is equal to the k-means cost functional for the set of rows of $V$,
where $V_{p, :}$ is the $p$th row of $V$ (and similarly for $V_{q, :}$).
Therefore, applying the k-means algorithm to the rows of $V$ will minimize an
upper bound on the largest principal angle between $\image{V}$ and
$\image{\pmat}$.

The advantage of using k-means compared to QR decomposition-based clustering is
in that the latter can only, given $V \in \bbR^{\vnum \times r_P}$, return a
partition with $r_P$ clusters.
On the other hand,
k-means clustering can return a partition with any number of clusters $\cnum$.
This makes it more efficient when $r_P \ll \cnum$ and projection-based MOR
method already generates a good subspace $\image{V}$.

For multi-agent systems~\eqref{eq:BGM20:general-mas} with agents of order $n$,
we have the matrices $V$ and $W$ as in~\eqref{eq:BGM20:VWPn}.
QR decomposition-based clustering can then be extended as in
Algorithm~2 from~\cite{morMliGB16}
by clustering the block-columns of $V\tran$ (or $W\tran$).
For the k-means algorithm, we can show in a similar way as in the
single-integrator case that clustering the block-rows leads to minimizing an
upper bound of the largest principal angle.
Therefore, k-means can be directly applied to the set of block-rows of $V$ or
$W$.

Note that the approach is not limited to the two clustering algorithms mentioned
here.
Any clustering algorithm over the set of (block-)rows of $V$ or $W$ can be used.
In particular, if only partitions with certain properties are wanted
(e.g., those that only cluster neighboring agents),
then special clustering algorithms could be used
(e.g., agglomerative clustering taking into account the connectivity of the
graph).

\section{Clustering for nonlinear multi-agent
  systems}\label{sec:BGM20:nonlinear-clustering}
In this section,
we extend the approach from the previous section to a class of nonlinear
multi-agent systems.
We describe the class of multi-agent systems in
Section~\ref{sec:BGM20:nonlin-cl-sys}.
Next,
in Section~\ref{sec:BGM20:nonlin-cl-proj},
we show that clustering by projection preserves structure for this class of
systems.

\subsection{Nonlinear multi-agent systems}\label{sec:BGM20:nonlin-cl-sys}
\begin{subequations}\label{eq:BGM20:nonlin-mas}
  Here,
  we consider a class of nonlinear multi-agent systems.
  In particular,
  let the dynamics of the $i$th agent,
  for $i = 1, 2, \dots, \vnum$,
  be defined by the control-affine system
  \begin{align}
    \dot{x}_i(t) & = A(x_i(t)) + B(x_i(t)) v_i(t), \\
    z_i(t) & = C(x_i(t)),
  \end{align}
  with functions
  $\fundef{A}{\Rn}{\Rn}$,
  $\fundef{B}{\Rn}{\Rnm}$,
  $\fundef{C}{\Rn}{\Rp}$,
  state $x_i(t) \in \Rn$,
  input $v_i(t) \in \Rm$, and
  output $z_i(t) \in \Rp$.
  Furthermore,
  let the interconnections be
  \begin{align}
    \mass_i v_i(t)
    & =
      \sum_{j = 1}^{\vnum}{\weight_{ij} K(z_i(t), z_j(t))}
      + \sum_{k = 1}^{\inum}{\iw_{ik} u_k(t)},
  \end{align}
  for $i = 1, 2, \dots, \vnum$,
  with inertias $\mass_i > 0$ and $\mmat = \diag{\mass_i}$,
  coupling $\fundef{K}{\Rp \times \Rp}{\Rm}$,
  external input $u_k(t) \in \Rm$, $k = 1, 2, \dots, \inum$,
  where
  $\amat = [\weight_{ij}]$ is the adjacency matrix of the graph $\graph$, and
  $\imat = [\iw_{ik}]$.
  Additionally, let the external output be
  \begin{align}
    y_{\ell}(t)
    & =
      \sum_{j = 1}^{\vnum}{\ow_{\ell j} z_j(t)},
  \end{align}
  with $\omat = [\ow_{\ell j}]$.
  We assume functions $A, B, C, K$ are continuous and that there is a unique
  global solution $x(t) = \col{x_1(t), x_2(t), \dots, x_{\vnum}(t)}$ for any
  admissible $u(t)$.
\end{subequations}

\subsection{Clustering by projection}\label{sec:BGM20:nonlin-cl-proj}
We want to find the form of the reduced order model obtained from Galerkin
projection with $V = \pmat \otimes I_n$.
We can rewrite~\eqref{eq:BGM20:nonlin-mas} to
\begin{align*}
  (\mmat \otimes I_n) \dot{x}(t)
  & = f(x(t), u(t)), \\
  y(t)
  & = g(x(t)),
\end{align*}
for some functions $f$ and $g$.
The reduced model is
\begin{align}
  \label{eq:BGM20:nonlin-cl-ode}
  \myparen*{\pmat\tran \mmat \pmat \otimes I_n} \dot{\hx}(t)
  & =
    \myparen*{\pmat\tran \otimes I_n}
    f\myparen*{(\pmat \otimes I_n) \hx(t), u(t)}, \\
  \nonumber
  \hy(t)
  & =
    g\myparen*{(\pmat \otimes I_n) \hx(t)},
\end{align}
with $\hx(t) = \col{\hx_1(t), \hx_2(t), \dots, \hx_{\cnum}(t)}$ and
$\hx_i(t) \in \Rn$.
Let $\pi(j) \in \{1, 2, \dots, \cnum\}$ be
such that $j \in \clu_{\pi(j)}$, for $j \in \{1, 2, \dots, \vnum\}$.
Premultiplying~\eqref{eq:BGM20:nonlin-cl-ode} with $e_{\imath}\tran \otimes I_n$
for some $\imath \in \{1, 2, \dots, \cnum\}$,
we find
\begin{align*}
  &
    \hat{\mass}_{\imath}
    \dot{\hx}_{\imath}(t) \\
  & =
    \sum_{i \in \clu_{\imath}}
    \myparen*{
      \mass_i A\myparen*{\hx_{\imath}(t)}
      + B\myparen*{\hx_{\imath}(t)}
      \myparen*{
        \sum_{j = 1}^{\vnum} \weight_{ij}
          K\myparen*{C\myparen*{\hx_{\imath}(t)}, C\myparen*{\hx_{\pi(j)}(t)}}
        + \sum_{k = 1}^{\inum}{\iw_{i k} u_{k}(t)}
      }
    }
    \displaybreak[0] \\
  & =
    \hat{\mass}_{\imath} A\myparen*{\hx_{\imath}(t)} \\
  & \qquad
    + B\myparen*{\hx_{\imath}(t)}
    \myparen*{
      \sum_{i \in \clu_{\imath}} \sum_{j = 1}^{\vnum}
        \weight_{ij}
        K\myparen*{C\myparen*{\hx_{\imath}(t)}, C\myparen*{\hx_{\pi(j)}(t)}}
      + \sum_{i \in \clu_{\imath}} \sum_{k = 1}^{\inum}
        \iw_{i k} u_{k}(t)
    }
    \displaybreak[0] \\
  & =
    \hat{\mass}_{\imath} A\myparen*{\hx_{\imath}(t)} \\
  & \qquad
    + B\myparen*{\hx_{\imath}(t)}
    \myparen*{
      \sum_{\jmath = 1}^{\cnum}
      \sum_{i \in \clu_{\imath}}
      \sum_{j \in \clu_{\jmath}}
        \weight_{ij}
        K\myparen*{C\myparen*{\hx_{\imath}(t)}, C\myparen*{\hx_{\jmath}(t)}}
      + \sum_{k = 1}^{\inum} \sum_{i \in \clu_{\imath}}
        \iw_{i k} u_{k}(t)
    } \\
  & =
    \hat{\mass}_{\imath} A\myparen*{\hx_{\imath}(t)}
    + B\myparen*{\hx_{\imath}(t)}
    \myparen*{
      \sum_{\jmath = 1}^{\cnum}
        \hat{\weight}_{\imath \jmath}
        K\myparen*{C\myparen*{\hx_{\imath}(t)}, C\myparen*{\hx_{\jmath}(t)}}
      + \sum_{k = 1}^{\inum}
        \hat{\iw}_{\imath k} u_{k}(t)
    },
\end{align*}
for
\begin{align*}
  \hat{\mass}_{\imath}
  =
    \sum_{i \in \clu_{\imath}} \mass_i, \quad
  \hat{\weight}_{\imath \jmath}
  & =
    \sum_{i \in \clu_{\imath}} \sum_{j \in \clu_{\jmath}} \weight_{ij}, \quad
  \hat{\iw}_{\imath k}
  =
    \sum_{i \in \clu_{\imath}} \iw_{i k}.
\end{align*}
Defining
$\hat{\mmat} \coloneqq \diag{\hat{\mass}_{\imath}}$,
$\hat{\amat} \coloneqq [\hat{\weight}_{\imath \jmath}]$, and
$\hat{\imat} \coloneqq [\hat{\iw}_{\imath k}]$,
we see that
$\hat{\mmat} = \pmat\tran \mmat \pmat$,
$\hat{\amat} = \pmat\tran \amat \pmat$, and
$\hat{\imat} = \pmat\tran \imat$.
For the output,
we have
\begin{align*}
  \hy_{\ell}(t)
  & =
    \sum_{j = 1}^{\vnum}{\ow_{\ell j} C\myparen*{\hx_{\pi(j)}(t)}}
    =
    \sum_{\jmath = 1}^{\cnum}
    \sum_{j \in \clu_{\jmath}}{\ow_{\ell j} C\myparen*{\hx_{\jmath}(t)}}
    =
    \sum_{\jmath = 1}^{\cnum}{
      \hat{\ow}_{\ell \jmath} C\myparen*{\hx_{\jmath}(t)}},
\end{align*}
where
\begin{align*}
  \hat{\ow}_{\ell \jmath}
  & =
    \sum_{j \in \clu_{\jmath}} \ow_{\ell j}.
\end{align*}
Thus, for $\hat{\omat} \coloneqq [\hat{\ow}_{\ell \jmath}]$, we have
$\hat{\omat} = \omat \pmat$.
Therefore, we showed how to construct a reduced model of the same structure as
the original multi-agent system.
Based on this, to find a good partition, we can apply any projection-based
MOR method for nonlinear systems (e.g., proper orthogonal
decomposition~\cite{morHinV05}) and cluster the block-rows of the matrix used to
project the system.

\section{Numerical examples}\label{sec:BGM20:example}
Here,
we demonstrate our approach for different network examples,
beginning with a small linear multi-agent system in
Section~\ref{sec:BGM20:example-small}.
Next,
in Section~\ref{sec:BGM20:example-vanderpol},
we use the Van der Pol oscillator network.

The source code of the implementations used to compute the presented results can
be obtained from
\begin{center}
  \url{https://doi.org/10.5281/zenodo.3924653}
\end{center}
and is authored by Petar~Mlinari\'c.

\subsection{Small network example}\label{sec:BGM20:example-small}
\begin{figure}[tb]
  \centering
  \begin{tikzpicture}[x=1em, y=1em, scale=5]
    \def\xyscale{1.5}
    \def\nodesize{1.7em}
    \def\lw{0.1em}
    \def\nodeweightsep{0.15em}

    \foreach \name/\pos in {%
      {1/(0, 0)}, {2/(0, 1)}, {3/(0, 2)}, {4/(0, 3)},
      {5/(\xyscale, 2)}, {6/(\xyscale, 1)},
      {7/(2 * \xyscale, 1)}, {8/(2 * \xyscale, 2)},
      {9/(3 * \xyscale, 1)}, {10/(3 * \xyscale, 0)}%
    }
      \node[
        shape=circle,
        draw,
        line width=\lw,
        inner sep=0em,
        outer sep=0em,
        minimum size=\nodesize,
      ] (\name) at \pos {$\name$};

    \begin{scope}[line width=\lw]
      \foreach \src/\dest/\weight in {%
        {1/6/5}, {2/6/2}, {3/4/1}, {3/5/2}, {4/5/5},
        {5/6/2}, {5/8/7}, {6/7/6},
        {7/8/1}, {7/9/1}, {7/10/1}%
      }
        \draw (\src) --
          node[
            shape=circle,
            fill=white,
            inner sep=\nodeweightsep,
            font=\footnotesize,
          ] {$\weight$} (\dest);
      \foreach \src/\dest/\weight in {{2/5/3}, {6/3/3}, {7/5/6}, {6/8/7}}
        \draw (\src) --
          node[
            near end,
            shape=circle,
            fill=white,
            inner sep=\nodeweightsep,
            font=\footnotesize
          ] {$\weight$} (\dest);
    \end{scope}
  \end{tikzpicture}
  \caption{Undirected weighted graph with $10$ vertices from~\cite{morMonTC14}}%
  \label{fig:small-network}
\end{figure}
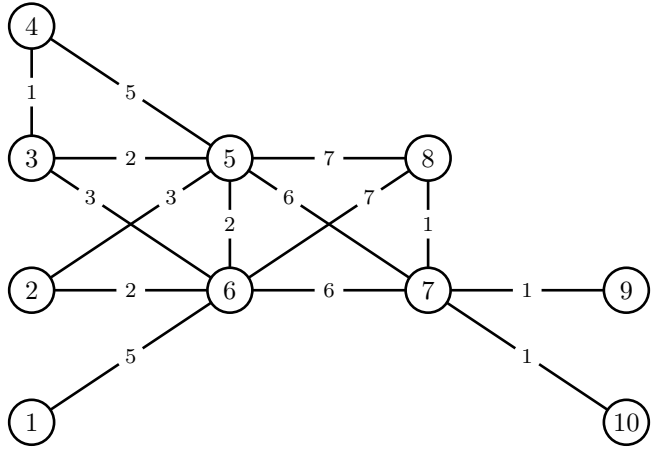

To illustrate distance to optimality,
we use the leader-follower multi-agents system example from~\cite{morMonTC14}
with $10$ single-integrator agents shown in Figure~\ref{fig:small-network},
where we can compute the $\Htwo$ and $\Hinf$ errors for all possible partitions.
The Laplacian and input matrices are
\begin{gather*}
	\lmat =
	\begin{bsmallmatrix*}[r]
		5 & 0 & 0 & 0 & 0 & -5 & 0 & 0 & 0 & 0 \\
		0 & 5 & 0 & 0 & -3 & -2 & 0 & 0 & 0 & 0 \\
		0 & 0 & 6 & -1 & -2 & -3 & 0 & 0 & 0 & 0 \\
		0 & 0 & -1 & 6 & -5 & 0 & 0 & 0 & 0 & 0 \\
		0 & -3 & -2 & -5 & 25 & -2 & -6 & -7 & 0 & 0 \\
		-5 & -2 & -3 & 0 & -2 & 25 & -6 & -7 & 0 & 0 \\
		0 & 0 & 0 & 0 & -6 & -6 & 15 & -1 & -1 & -1 \\
		0 & 0 & 0 & 0 & -7 & -7 & -1 & 15 & 0 & 0 \\
		0 & 0 & 0 & 0 & 0 & 0 & -1 & 0 & 1 & 0 \\
		0 & 0 & 0 & 0 & 0 & 0 & -1 & 0 & 0 & 1
	\end{bsmallmatrix*},
	\quad
	\imat =
	\begin{bsmallmatrix}
		0 & 0 \\
		0 & 0 \\
		0 & 0 \\
		0 & 0 \\
		0 & 0 \\
		1 & 0 \\
		0 & 1 \\
		0 & 0 \\
		0 & 0 \\
		0 & 0
	\end{bsmallmatrix},
\end{gather*}
and we chose the edge ordering and orientation such that the incidence and
edge-weights matrices are
\begin{gather*}
  \setcounter{MaxMatrixCols}{20}
	\rmat =
	\begin{bsmallmatrix*}[r]
		-1 & 0 & 0 & 0 & 0 & 0 & 0 & 0 & 0 & 0 & 0 & 0 & 0 & 0 & 0 \\
		0 & -1 & -1 & 0 & 0 & 0 & 0 & 0 & 0 & 0 & 0 & 0 & 0 & 0 & 0 \\
		0 & 0 & 0 & -1 & -1 & -1 & 0 & 0 & 0 & 0 & 0 & 0 & 0 & 0 & 0 \\
		0 & 0 & 0 & 1 & 0 & 0 & -1 & 0 & 0 & 0 & 0 & 0 & 0 & 0 & 0 \\
		0 & 1 & 0 & 0 & 1 & 0 & 1 & -1 & -1 & -1 & 0 & 0 & 0 & 0 & 0 \\
		1 & 0 & 1 & 0 & 0 & 1 & 0 & 1 & 0 & 0 & -1 & -1 & 0 & 0 & 0 \\
		0 & 0 & 0 & 0 & 0 & 0 & 0 & 0 & 1 & 0 & 1 & 0 & -1 & -1 & -1 \\
		0 & 0 & 0 & 0 & 0 & 0 & 0 & 0 & 0 & 1 & 0 & 1 & 1 & 0 & 0 \\
		0 & 0 & 0 & 0 & 0 & 0 & 0 & 0 & 0 & 0 & 0 & 0 & 0 & 1 & 0 \\
		0 & 0 & 0 & 0 & 0 & 0 & 0 & 0 & 0 & 0 & 0 & 0 & 0 & 0 & 1
	\end{bsmallmatrix*}
\end{gather*}
and $\wmat = \diag{5, 3, 2, 1, 2, 3, 5, 2, 6, 7, 6, 7, 1, 1, 1}$,
respectively.
The output matrix is $\omat = \wmat^{1/2} \rmat\tran$.

For this example, we focus on partitions with five clusters.
There are in total $42\,525$ such partitions.
Table~\ref{tab:BGM20:top} shows the $15$ best partitions with respect to the
$\Htwo$ and $\Hinf$ errors.
\begin{table}[tb]
	\caption{Top $15$ partitions with $5$ clusters by $\Htwo$ error and
    $\Hinf$ error for reducing the multi-agent system in
    Section~\ref{sec:BGM20:example-small}}\label{tab:BGM20:top}
  \centering
  \rowcolors{1}{}{black!10}
  \begin{tabular}{ccc}
		Rank & Relative $\Htwo$ error & Partition \\
		\midrule
    $ 1$ & $0.128053$ & $\{\{1, 8\}, \{2, 3, 4, 9, 10\}, \{5\}, \{6\}, \{7\}\}$ \\
    $ 2$ & $0.131311$ & $\{\{1, 2, 3, 4\}, \{5, 8\}, \{6\}, \{7\}, \{9, 10\}\}$ \\
    $ 3$ & $0.137466$ & $\{\{1, 2, 3, 4, 9, 10\}, \{5\}, \{6\}, \{7\}, \{8\}\}$ \\
    $ 4$ & $0.137473$ & $\{\{1, 3, 8\}, \{2, 4, 9, 10\}, \{5\}, \{6\}, \{7\}\}$ \\
    $ 5$ & $0.143700$ & $\{\{1, 5, 8\}, \{2, 3, 4\}, \{6\}, \{7\}, \{9, 10\}\}$ \\
    $ 6$ & $0.145900$ & $\{\{1, 2, 3\}, \{4, 9, 10\}, \{5, 8\}, \{6\}, \{7\}\}$ \\
    $ 7$ & $0.146196$ & $\{\{1, 8\}, \{2, 3, 4, 9\}, \{5, 10\}, \{6\}, \{7\}\}$ \\
    $ 8$ & $0.146196$ & $\{\{1, 8\}, \{2, 3, 4, 10\}, \{5, 9\}, \{6\}, \{7\}\}$ \\
    $ 9$ & $0.147022$ & $\{\{1, 2, 3, 8\}, \{4, 9, 10\}, \{5\}, \{6\}, \{7\}\}$ \\
    $10$ & $0.149240$ & $\{\{1, 8, 10\}, \{2, 3, 4, 9\}, \{5\}, \{6\}, \{7\}\}$ \\
    $11$ & $0.149240$ & $\{\{1, 8, 9\}, \{2, 3, 4, 10\}, \{5\}, \{6\}, \{7\}\}$ \\
    $12$ & $0.149654$ & $\{\{1, 8\}, \{2, 4, 9, 10\}, \{3, 5\}, \{6\}, \{7\}\}$ \\
    $13$ & $0.150440$ & $\{\{1, 5\}, \{2, 3, 4, 9, 10\}, \{6\}, \{7\}, \{8\}\}$ \\
    $14$ & $0.150654$ & $\{\{1, 3\}, \{2, 4, 9, 10\}, \{5, 8\}, \{6\}, \{7\}\}$ \\
    $15$ & $0.151684$ & $\{\{1, 2, 8\}, \{3, 4, 9, 10\}, \{5\}, \{6\}, \{7\}\}$ \\
  \end{tabular}

  \vspace*{1em}

  \rowcolors{1}{}{black!10}
  \begin{tabular}{ccc}
		Rank & Relative $\Hinf$ error & Partition \\
		\midrule
    $ 1$ & $0.253975$ & $\{\{1, 3, 5, 8\}, \{2, 4\}, \{6\}, \{7\}, \{9, 10\}\}$ \\
    $ 2$ & $0.254376$ & $\{\{1, 2, 5, 8\}, \{3, 4\}, \{6\}, \{7\}, \{9, 10\}\}$ \\
    $ 3$ & $0.254818$ & $\{\{1, 5, 8\}, \{2, 3, 4\}, \{6\}, \{7\}, \{9, 10\}\}$ \\
    $ 4$ & $0.259483$ & $\{\{1, 2, 3, 5, 8\}, \{4\}, \{6\}, \{7\}, \{9, 10\}\}$ \\
    $ 5$ & $0.260859$ & $\{\{1, 2, 4\}, \{3, 5, 8\}, \{6\}, \{7\}, \{9, 10\}\}$ \\
    $ 6$ & $0.262244$ & $\{\{1, 2, 3, 4\}, \{5, 8\}, \{6\}, \{7\}, \{9, 10\}\}$ \\
    $ 7$ & $0.266387$ & $\{\{1, 3, 4\}, \{2, 5, 8\}, \{6\}, \{7\}, \{9, 10\}\}$ \\
    $ 8$ & $0.273663$ & $\{\{1, 4\}, \{2, 3, 5, 8\}, \{6\}, \{7\}, \{9, 10\}\}$ \\
    $ 9$ & $0.276919$ & $\{\{1, 4, 5, 8\}, \{2, 3\}, \{6\}, \{7\}, \{9, 10\}\}$ \\
    $10$ & $0.286961$ & $\{\{1, 3, 4, 5, 8\}, \{2\}, \{6\}, \{7\}, \{9, 10\}\}$ \\
    $11$ & $0.288414$ & $\{\{1, 2, 3\}, \{4, 5, 8\}, \{6\}, \{7\}, \{9, 10\}\}$ \\
    $12$ & $0.293773$ & $\{\{1, 5\}, \{2, 3, 4, 8\}, \{6\}, \{7\}, \{9, 10\}\}$ \\
    $13$ & $0.294028$ & $\{\{1, 2, 3, 4, 8\}, \{5\}, \{6\}, \{7\}, \{9, 10\}\}$ \\
    $14$ & $0.299845$ & $\{\{1, 2\}, \{3, 4, 5, 8\}, \{6\}, \{7\}, \{9, 10\}\}$ \\
    $15$ & $0.305583$ & $\{\{1, 2, 4, 8\}, \{3, 5\}, \{6\}, \{7\}, \{9, 10\}\}$
  \end{tabular}
\end{table}

First, we used IRKA to find a reduced model of order $\cnum = 5$.
It found a reduced model with relative $\Htwo$ error of
$3.30412 \times 10^{-2}$, which is $3.88$ times better than the best partition.
The partition resulting from QR decomposition-based clustering applied to IRKA's
$V$ matrix is
\begin{equation*}
	\{\{1, 3\}, \{2, 4, 9, 10\}, \{5, 8\}, \{6\}, \{7\}\},
\end{equation*}
with the associated relative $\Htwo$ error of $0.150654$.
It is more than $4$ times worse than using IRKA,
but note that this partition is the $14$th best partition and
that the best partition produces about $1.18$ times better error.
Using k-means clustering gives
\begin{equation*}
  \{\{1, 2, 3\}, \{4, 9, 10\}, \{5, 8\}, \{6\}, \{7\}\},
\end{equation*}
with relative $\Htwo$ error of $0.1459$ and taking the $6$th place.

We notice by~\eqref{eq:BGM20:VWP} that $W$ can also be used to find a good
partition.
In this example, QR decomposition-based clustering returns the partition
\begin{equation*}
	\{\{1, 2, 3, 9, 10\}, \{4, 8\}, \{5\}, \{6\}, \{7\}\},
\end{equation*}
with the relative $\Htwo$ error $0.179746$,
which is worse than using only $V$ from IRKA\@.
Using k-means clustering returns
\begin{equation*}
  \{\{1, 2, 3, 4, 8\}, \{5\}, \{6\}, \{7\}, \{9, 10\}\}
\end{equation*}
with the relative $\Htwo$ error $0.156788$,

Using the first five left singular vectors of
$
\begin{bmatrix}
  V & W
\end{bmatrix}
$
to take into account both $V$ and $W$,
using QR decomposition-based clustering produces
\begin{equation*}
  \{\{1, 2, 3, 4, 8\}, \{5\}, \{6\}, \{7\}, \{9, 10\}\}
\end{equation*}
with the relative $\Htwo$ error $0.189487$,
which further increases the error.
On the other hand,
k-means clustering gives us
\begin{equation*}
  \{\{1, 2, 3, 4\}, \{5, 8\}, \{6\}, \{7\}, \{9, 10\}\},
\end{equation*}
which is the second best partition in terms of the $\Htwo$ error and sixth best
in terms of the $\Hinf$ error.
Furthermore,
using balanced truncation instead of IRKA produces the same partition,
using either of the two clustering algorithms and the three choices of matrices.

Therefore, at least in this example, clustering the rows of $V$ and/or $W$ gives
close to optimal partitions.
Additionally, k-means clustering performs better than QR decomposition-based
clustering.

\subsection{Van der Pol oscillators}\label{sec:BGM20:example-vanderpol}
\begin{subequations}\label{eq:BGM20:vanderpol}
  Here, we use the Van der Pol oscillator network example from~\cite{MasS18},
  where the agents are given by
  \begin{align}
    \dot{x}_{i, 1}(t)
    & =
      x_{i, 2}(t) + \sigma v_i(t), \\
    \dot{x}_{i, 2}(t)
    & =
      \mu \myparen*{1 - {x_{i, 1}(t)}^2} x_{i, 2}(t) - x_{i, 1}(t) - c v_i(t),
  \end{align}
  and interconnections by
  \begin{align}
    v_i(t)
    & =
      \sum_{j = 1}^{\vnum}{\weight_{ij}
      \myparen*{\myparen*{x_{i, 1}(t) - x_{j, 1}(t)}
      + \myparen*{x_{i, 2}(t) - x_{j, 2}(t)}}}
      + \sum_{k = 1}^{\inum}{\iw_{ik} u_k(t)},
  \end{align}
  with $\mu = 0.5$ and $\sigma = 0.1$. Additionally, we chose a larger $10
  \times 10$ grid graph ($\vnum = 100$), set the input matrix to be $\imat =
  e_1$ (i.e., one of the corner agents receives external input) and used $c =
  100$ to have synchronization.
\end{subequations}

Figure~\ref{fig:BGM20:vanderpol-train} shows the state trajectory of the system
for zero initial condition and input $u(t) = e^{-t}$,
using an adaptive BDF integrator producing $987$ snapshots.
We used these snapshots to find the POD modes, with associated singular values
shown in Figure~\ref{fig:BGM20:vanderpol-svdvals}.
\begin{figure}[tb]
  \centering
  \begin{tikzpicture}
    \begin{axis}[
        width=0.8\linewidth,
        height=0.5\linewidth,
        xlabel={Time},
        grid=major,
        legend entries={$x_{i, 1}$, $x_{i, 2}$},
        legend pos=north west,
      ]
      \addplot [%
          very thick,
          tab_blue,
        ]
        graphics [
          xmin=-0.05, xmax=20.05,
          ymin=-1.99629, ymax=2.00125,
        ] {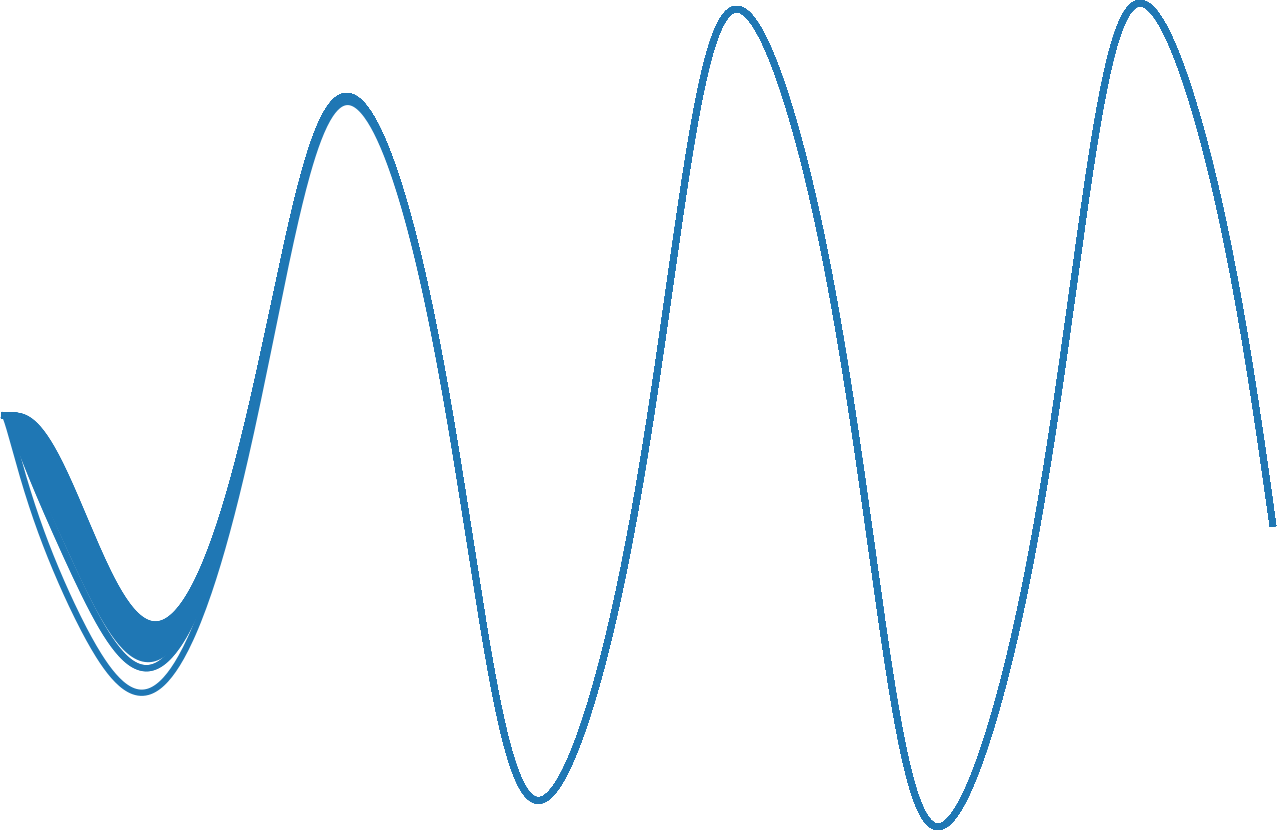};
      \addplot [%
          very thick,
          tab_orange,
        ]
        graphics [
          xmin=-0.05, xmax=20.05,
          ymin=-2.21841, ymax=2.22445,
        ] {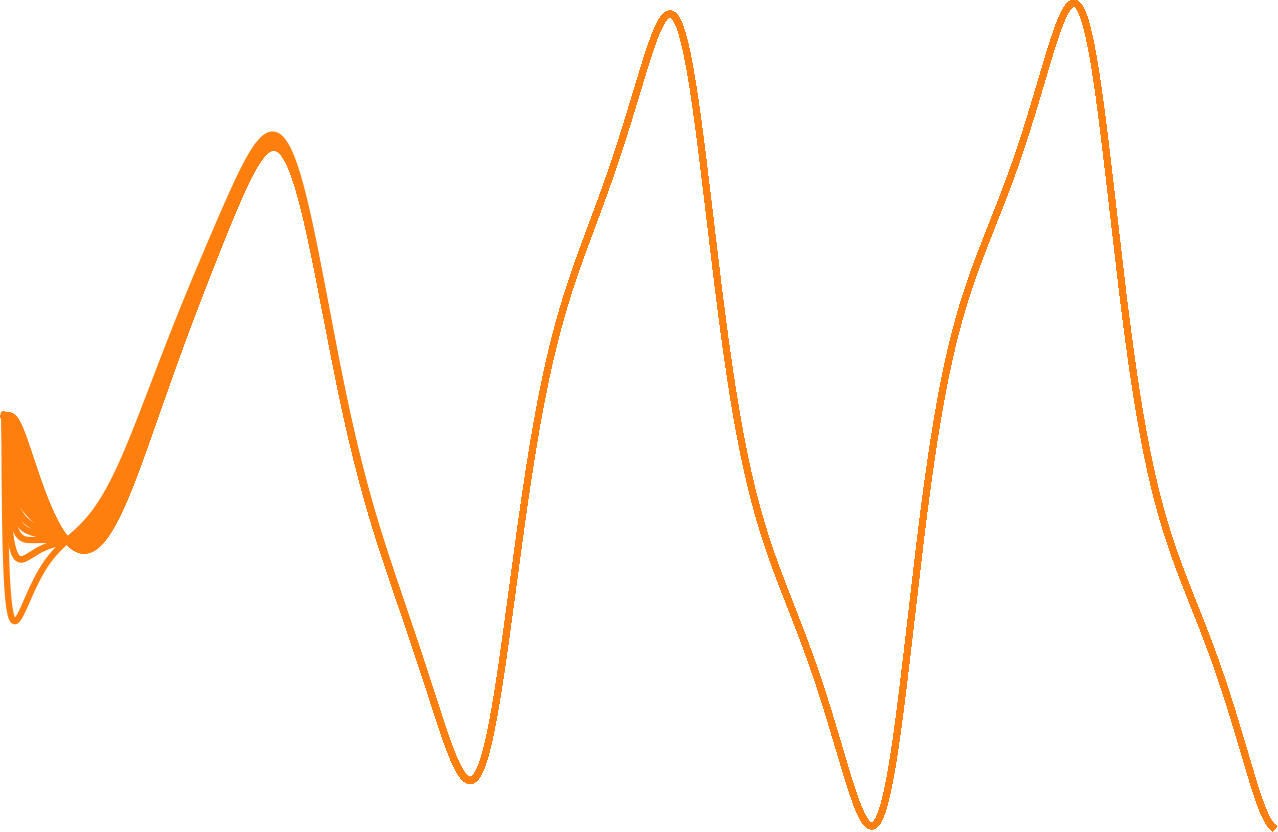};
    \end{axis}
  \end{tikzpicture}
  \caption{State trajectory of the Van der Pol oscillator
    network~\eqref{eq:BGM20:vanderpol} for zero initial condition and input
    $u(t) = e^{-t}$}%
  \label{fig:BGM20:vanderpol-train}
\end{figure}
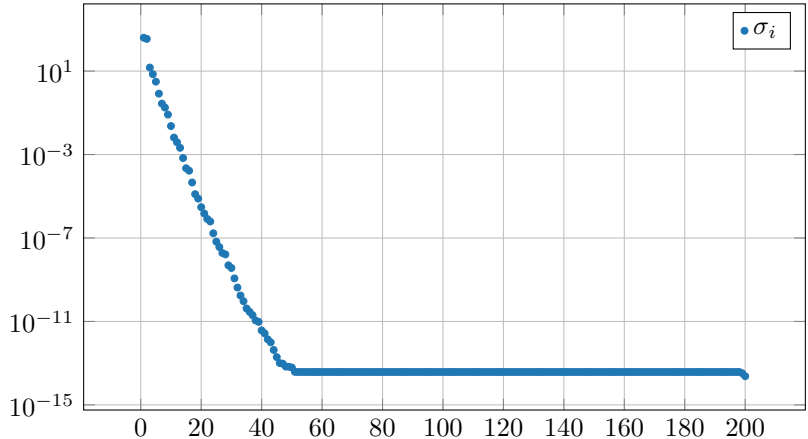
\begin{figure}[tb]
  \centering
  \begin{tikzpicture}
    \begin{semilogyaxis}[
        width=0.8\linewidth,
        height=0.5\linewidth,
        grid=major,
        legend entries={$\sigma_i$}
      ]
      \addplot [only marks, mark size=0.3ex, tab_blue]
        table [x index=0, y index=1]
        {code/fig/vanderpol-svdvals.dat};
    \end{semilogyaxis}
  \end{tikzpicture}
  \caption{POD singular values based on snapshots from
    Figure~\ref{fig:BGM20:vanderpol-train}}%
  \label{fig:BGM20:vanderpol-svdvals}
\end{figure}

\begin{figure}[tb]
  \centering
  \begin{tikzpicture}
    \begin{axis}[
        width=0.8\linewidth,
        height=0.5\linewidth,
        xlabel={Time},
        grid=major,
        legend entries={$x_{i, 1}$, $x_{i, 2}$},
        legend pos=north west,
      ]
      \addplot [%
          very thick,
          tab_blue,
        ]
        graphics [
          xmin=-0.05, xmax=20.05,
          ymin=-2.69061, ymax=2.80036,
        ] {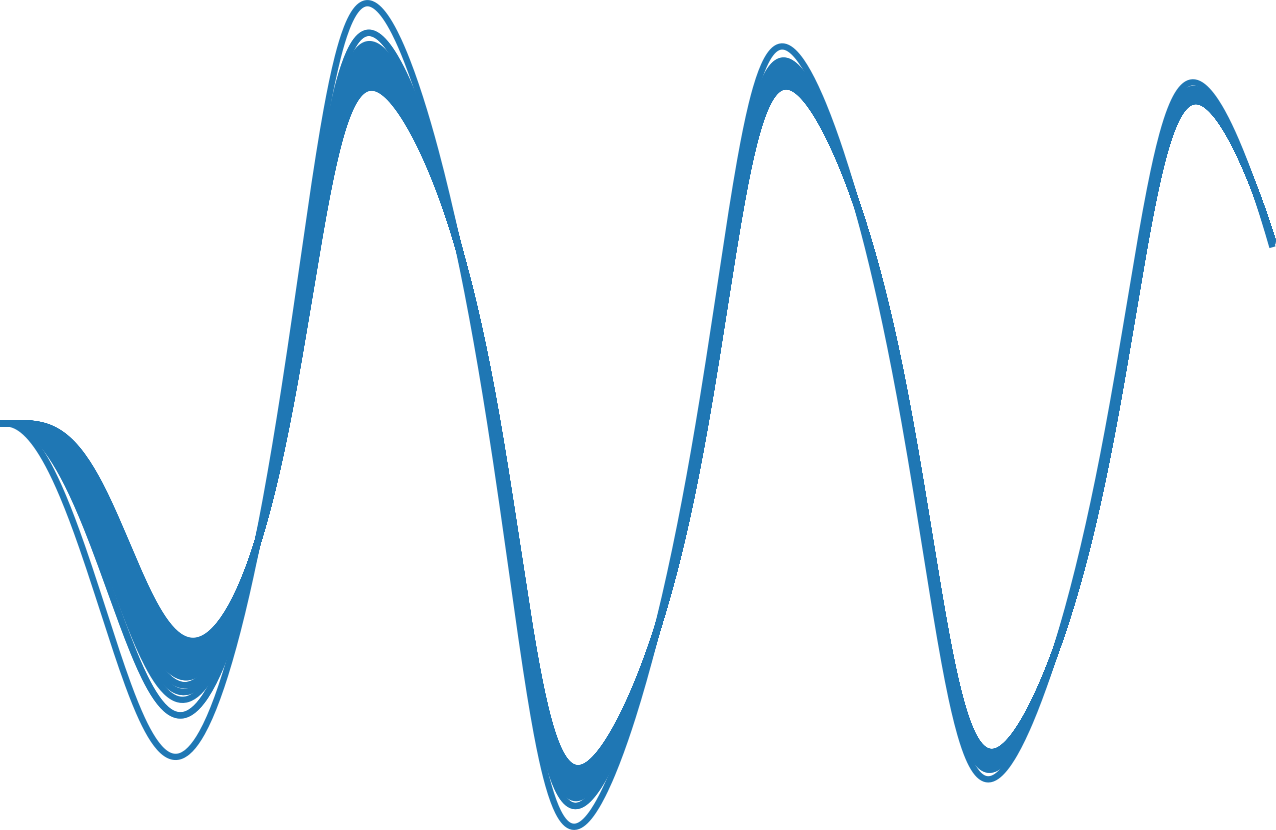};
      \addplot [%
          very thick,
          tab_orange,
        ]
        graphics [
          xmin=-0.05, xmax=20.05,
          ymin=-3.00636, ymax=2.94303,
        ] {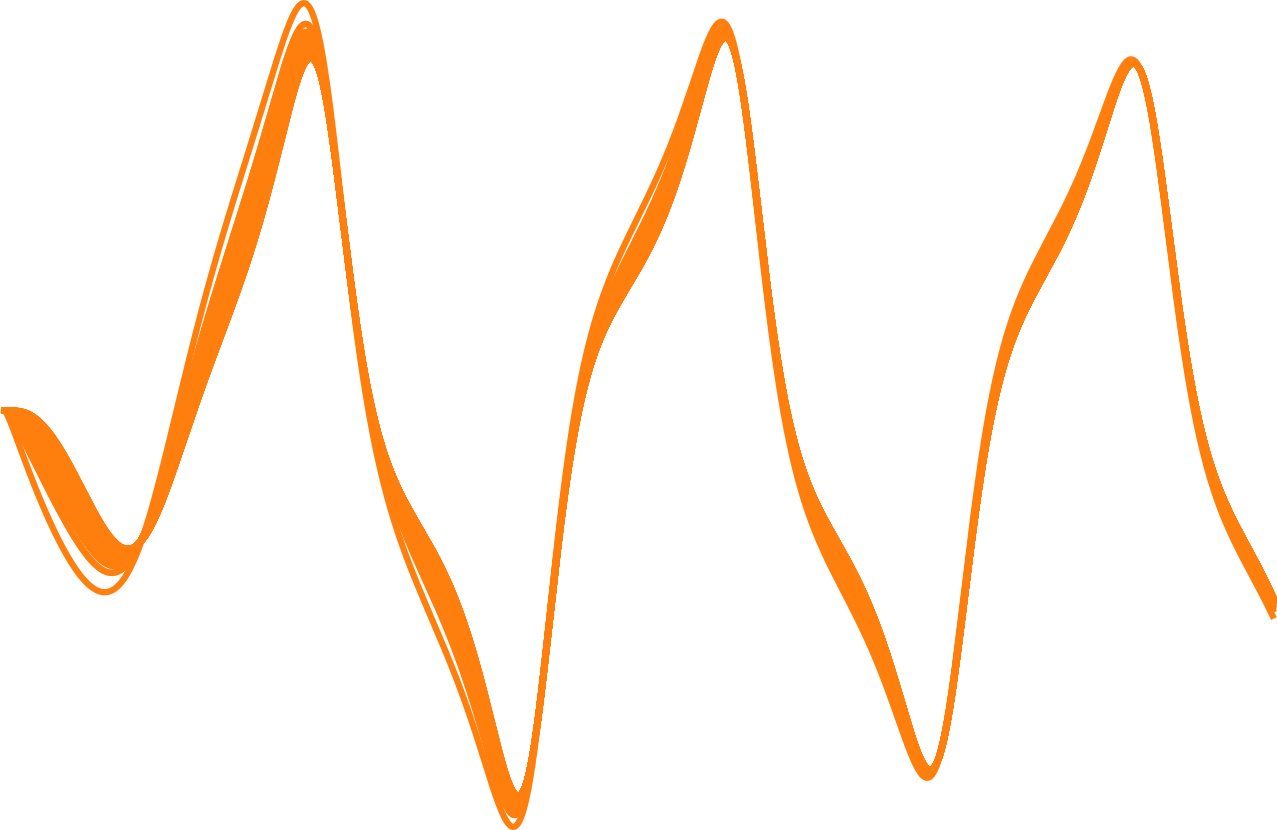};
    \end{axis}
  \end{tikzpicture}
  \caption{Van der Pol oscillator state trajectory for zero initial condition
    and input $u(t) = e^{-t/10} \sin{t}$}%
  \label{fig:BGM20:vanderpol-test}
\end{figure}

\begin{figure}[tb]
  \centering
  \begin{tikzpicture}
    \begin{axis}[
        width=0.8\linewidth,
        height=0.5\linewidth,
        xlabel={Time},
        grid=major,
        legend entries={%
          $x_{i, 1} - \pmat \hx_{i, 1}$,
          $x_{i, 2} - \pmat \hx_{i, 2}$},
        legend pos=north east,
        ymin=-0.105,
        ymax=0.105,
        yticklabel={$\tick$},
      ]
      \addplot [%
          very thick,
          tab_blue,
        ]
        graphics [
          xmin=-0.05, xmax=20.05,
          ymin=-0.10104, ymax=0.0682957,
        ] {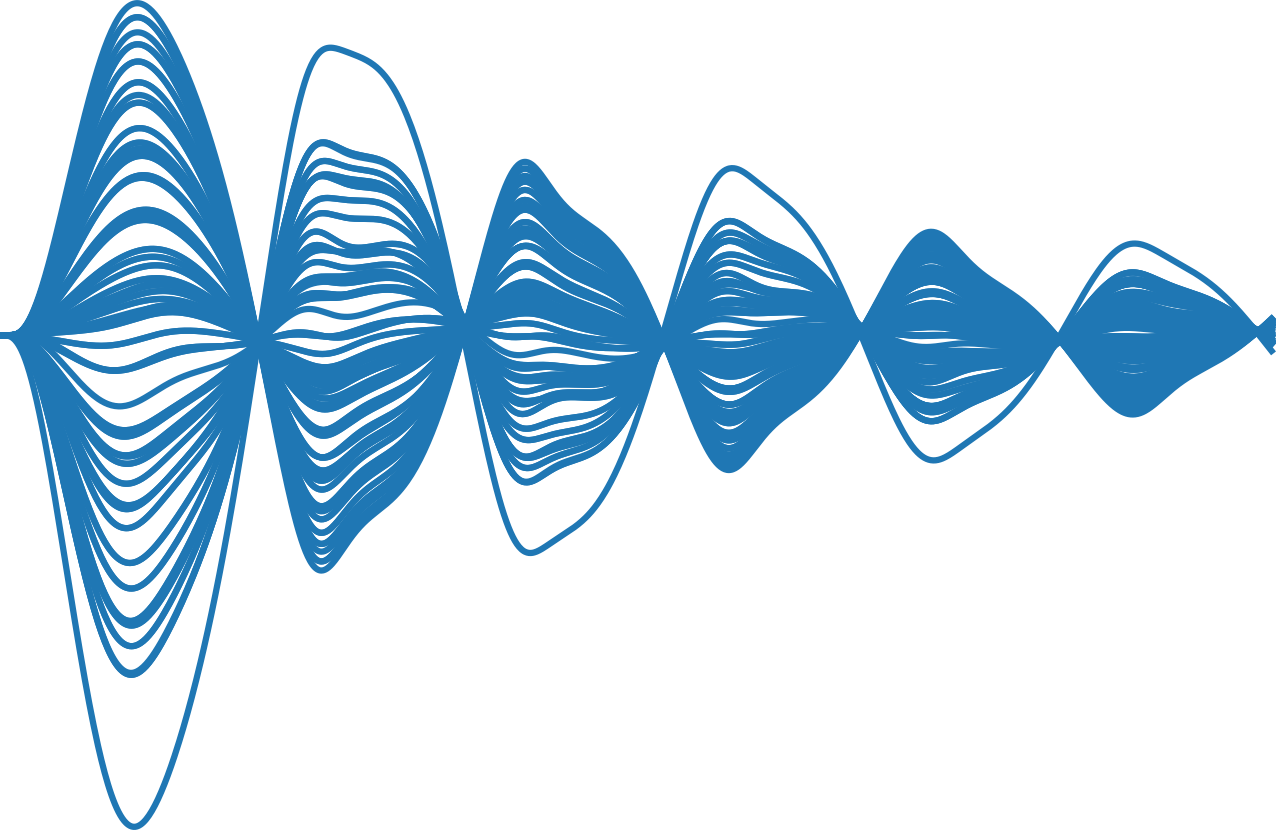};
      \addplot [%
          very thick,
          tab_orange,
        ]
        graphics [
          xmin=-0.05, xmax=20.05,
          ymin=-0.0813927, ymax=0.0861159,
        ] {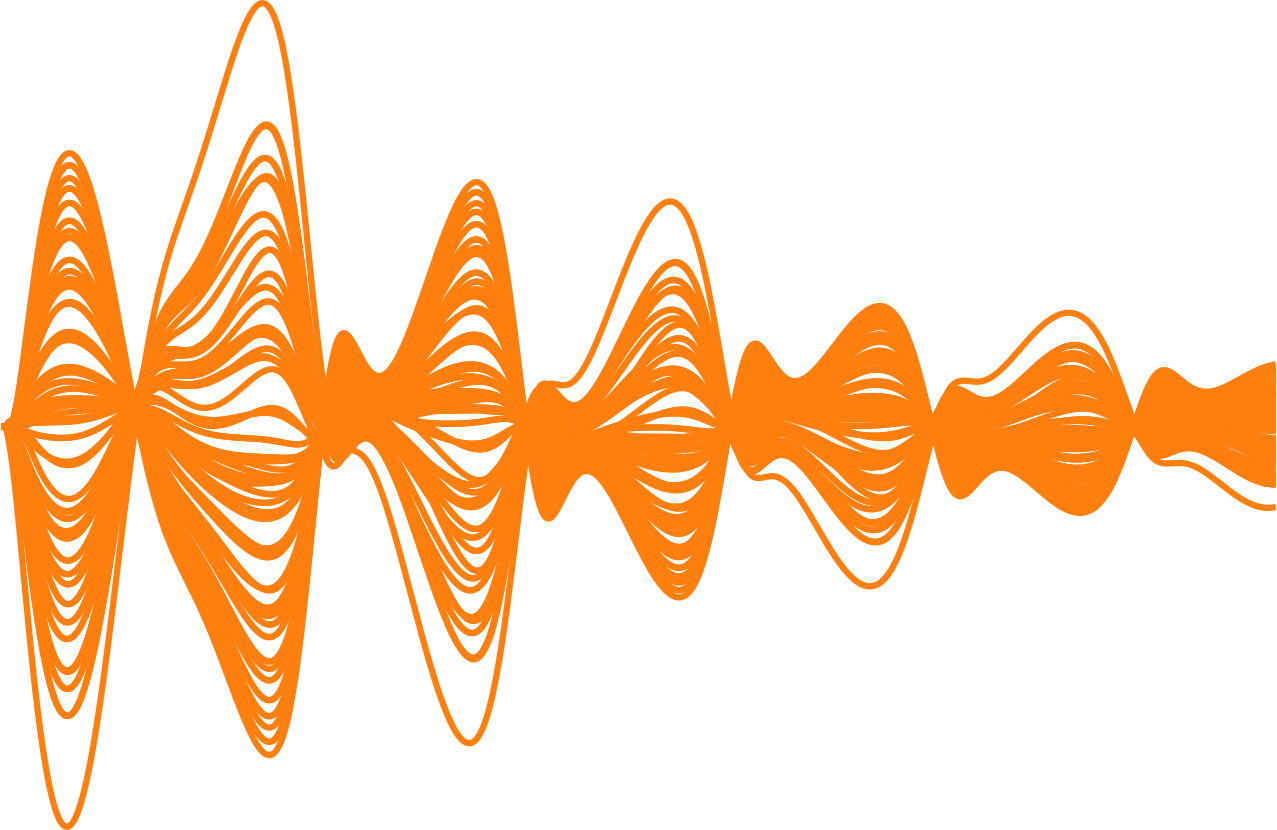};
    \end{axis}
  \end{tikzpicture}
  \caption{Van der Pol oscillator error when using k-means with the first two
    POD modes for zero initial condition and input $u(t) = e^{-t/10} \sin{t}$}%
  \label{fig:BGM20:vanderpol-err}
\end{figure}

\begin{figure}[tb]
  \centering
  \begin{tikzpicture}
    \begin{semilogyaxis}[
        width=0.8\linewidth,
        height=0.5\linewidth,
        xlabel={Reduced order},
        ylabel={Relative $\Ltwo$ error},
        grid=major,
      ]
      \addplot [very thick, mark=*, mark size=0.3ex, tab_blue]
        table [x index=0, y index=1]
        {code/fig/vanderpol-kmeans-error.dat};
    \end{semilogyaxis}
  \end{tikzpicture}
  \caption{Relative $\Ltwo$ error for zero initial condition and test input
    $u(t) = e^{-t/10} \sin{t}$ for k-means clustering using the first two POD
    modes}%
  \label{fig:BGM20:vanderpol-kmeans-error}
\end{figure}
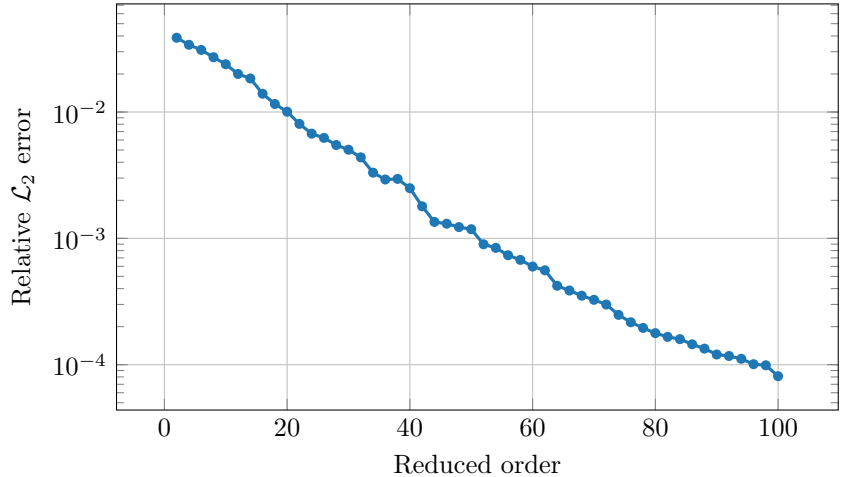

Changing the input to $u(t) = e^{-t/10} \sin{t}$ gives the trajectory in
Figure~\ref{fig:BGM20:vanderpol-test}.
Applying k-means clustering to the first two POD modes to generate $10$ clusters
produces a reduced model with the error trajectory in
Figure~\ref{fig:BGM20:vanderpol-err}.
We computed the relative $\Ltwo$ error for k-means clustering
using the first two POD modes with different number of clustering,
which can be seen in Figure~\ref{fig:BGM20:vanderpol-kmeans-error}.
For this example, we see that the error decays exponentially with the order of
the reduced model.

\section{Conclusions}\label{sec:BGM20:con}
We extended clustering by projection to a class of nonlinear multi-agent
systems and presented our clustering-based MOR method,
combining any projection-based MOR method and a clustering algorithm,
for reduction of multi-agent systems using graph partitions.
In particular, we motivated the use of the k-means algorithm.

Our numerical test for a small network shows that our algorithm finds close to
optimal partitions.
We demonstrated our method on a larger nonlinear oscillator network.

\bibliographystyle{plain}
\bibliography{refs}
\end{document}